\documentclass[12pt]{article}
\usepackage{amsthm,amssymb,amsmath}
\usepackage{epsfig}
\usepackage{verbatim}
\usepackage{fullpage}

\title{Planar digraphs of digirth five are 2-colorable}

\author{Ararat Harutyunyan\thanks{Research supported by a Digiteo postdoctoral scholarship.}\\
  {Laboratoire de Recherche en Informatique}\\
  {Universit\'{e} Paris-Sud}\\
  {91405 Orsay, France} \\
  email: {\tt ararat@lri.fr}
\and
  Bojan Mohar\thanks{Supported in part by an NSERC Discovery Grant (Canada),
  by the Canada Research Chair program, and by the
  Research Grant P1--0297 of ARRS (Slovenia).}~\thanks{On leave from:
  IMFM \& FMF, Department of Mathematics, University of Ljubljana, Ljubljana,
  Slovenia.}\\
  {Department of Mathematics}\\
  {Simon Fraser University}\\
  {Burnaby, B.C. V5A 1S6} \\
  email: {\tt mohar@sfu.ca}
}

\newtheorem{theorem}{Theorem}[section]
\newtheorem{lemma}[theorem]{Lemma}

\newtheorem{conjecture}[theorem]{Conjecture}
\newtheorem{claim}{Claim}

\newcommand{\DEF}[1]{{\em #1\/}}

\newcommand{\C}{{\cal C}}

\begin{document}

\maketitle

\begin{abstract}
Neumann-Lara (1985) and \v{S}krekovski conjectured that every
planar digraph with digirth at least three is 2-colorable. We
prove a relaxed version of this conjecture: every planar digraph
of digirth at least five is 2-colorable. The result also holds in
the setting of list colorings.
\end{abstract}


\section{Introduction}

Let $D$ be a digraph without cycles of length $\leq 2$, and let
$G$ be the underlying undirected graph of $D$. A function $f
\colon V(D) \to \{1,\dots,k \}$ is a \DEF{k-coloring} of the
digraph $D$ if $V_i = f^{-1}(i)$ is acyclic in $D$ for every
$i=1,\dots,k$. Here we treat the vertex set $V_i$ \DEF{acyclic} if
the induced subdigraph $D[V_i]$ contains no directed cycles (but
$G[V_i]$ may contain cycles). We say that $D$ is \DEF{$k$-colorable}
if it admits a $k$-coloring. The minimum $k$ for which $D$ is
$k$-colorable is called the \DEF{chromatic number} of $D$, and is
denoted by $\chi(D)$ (see Neumann-Lara \cite{N1982}).

The following conjecture was proposed independently by
Neumann-Lara \cite{N} and \v{S}krekovski (see~\cite{BFJKM2004}).

\begin{conjecture}
\label{conj:planar}
Every planar digraph $D$ with no directed cycles of length at most
$2$ is $2$-colorable.
\end{conjecture}

The \DEF{digirth} of a digraph is the length of its shortest directed
cycle ($\infty$ if $D$ is acyclic).
It is an easy consequence of 5-degeneracy  of planar graphs that
every planar digraph $D$ with digirth at least 3 has chromatic number at most~3.

There seem to be lack of methods to attack Conjecture
\ref{conj:planar}, and no nontrivial partial results are known.
The main result of this paper is the following theorem whose proof
is based on elaborate use of (nowadays standard) discharging
technique.

\begin{theorem}
\label{thm:main}
Every planar digraph that has digirth at least five is 2-colorable.
\end{theorem}

The proof of Theorem \ref{thm:main} is deferred until Section 4.
Actually, we shall prove an extended version in the setting of
list-colorings which we define next.

Let $\C$ be a finite set of colors. Given a digraph $D$, let $L:
v\mapsto L(v)\subseteq \C$ be a \DEF{list-assignment} for $D$,
which assigns to each vertex $v\in V(D)$ a set of colors. The set
$L(v)$ is called the \DEF{list} (or the set of \DEF{admissible
colors}) for $v$. We say $D$ is \DEF{$L$-colorable} if there is an
\DEF{$L$-coloring} of $D$, i.e., each vertex $v$ is assigned a
color from $L(v)$ such that every color class induces an acyclic
set in $D$. A \DEF{$k$-list-assignment} for $D$ is a
list-assignment $L$ such that $|L(v)|=k$ for every $v \in V(D)$.
We say that $D$ is \DEF{$k$-choosable} if it is $L$-colorable for
every $k$-list-assignment $L$.

\begin{theorem}
\label{thm:main choosable}
Every planar digraph of digirth at least five is 2-choosable.
\end{theorem}

The rest of the paper is devoted to the proof of Theorem \ref{thm:main choosable}.

\section{Unavoidable configurations}

In this section we provide a list of unavoidable configurations used in the proof of Theorem \ref{thm:main choosable}.
Orientations of edges are not important at this point, so we shall consider only undirected graphs throughout the whole section.

We define a \DEF{configuration} as a plane graph $C$ together with
a function $\delta: U \to \mathbb{N}$, where $U\subseteq V(C)$, such that
$\delta(v) \geq deg_C(v)$ for every $v \in V(C)$. A plane graph
$G$ \DEF{contains} the configuration $(C, U, \delta)$ if there is a
mapping $h: V(C) \to V(G)$ with the following properties:

\begin{enumerate}
\item[(i)] For every edge $ab \in E(C)$, $h(a)h(b)$ is an edge of $G$.
\item[(ii)] For every facial walk $a_1 \dots a_k$ in $C$, except for the unbounded face, the image $h(a_1) \dots h(a_k)$ is a facial walk in $G$.
\item[(iii)] For every $a \in U$, the degree of $h(a)$ in $G$ is equal to $\delta(a)$.
\item[(iv)] $h$ is \DEF{locally one-to-one}, i.e., it is one-to-one on the neighbors of each vertex of $V(C)$.
\end{enumerate}

\begin{figure}[t!]
   \centering
   \includegraphics[width=14cm]{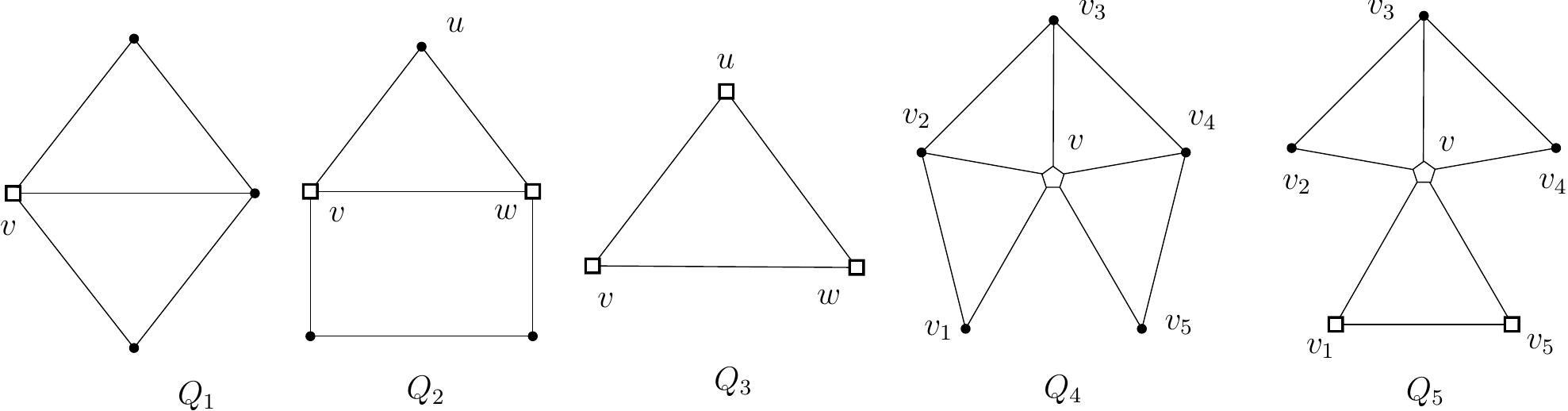}
   \caption{Configurations $Q_1$ to $Q_5$}
   \label{fig:1}
\end{figure}

\begin{figure}[t!]
   \centering
   \includegraphics[width=13.5cm]{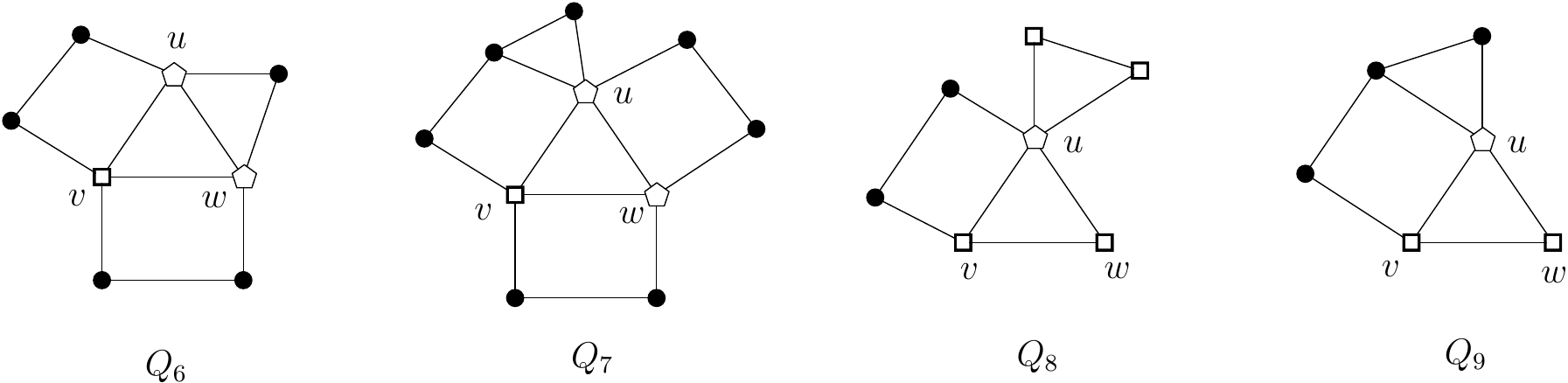}
   \caption{Configurations $Q_6$ to $Q_9$}
   \label{fig:2}
\end{figure}

\begin{figure}[t!]
   \centering
   \includegraphics[width=13.8cm]{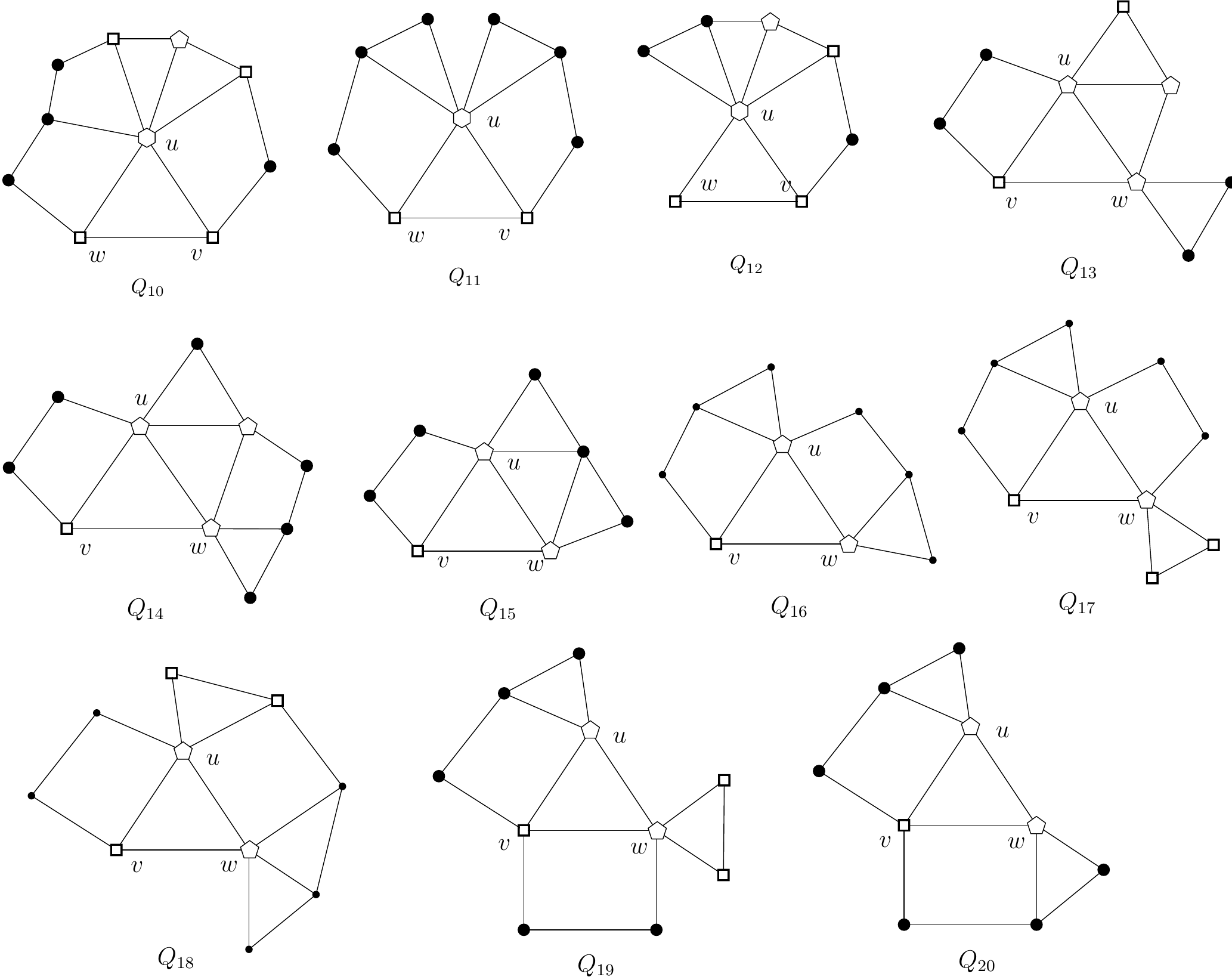}
   \caption{Configurations $Q_{10}$ to $Q_{20}$}
   \label{fig:3}
\end{figure}

\begin{figure}[t!]
   \centering
   \includegraphics[width=12.8cm]{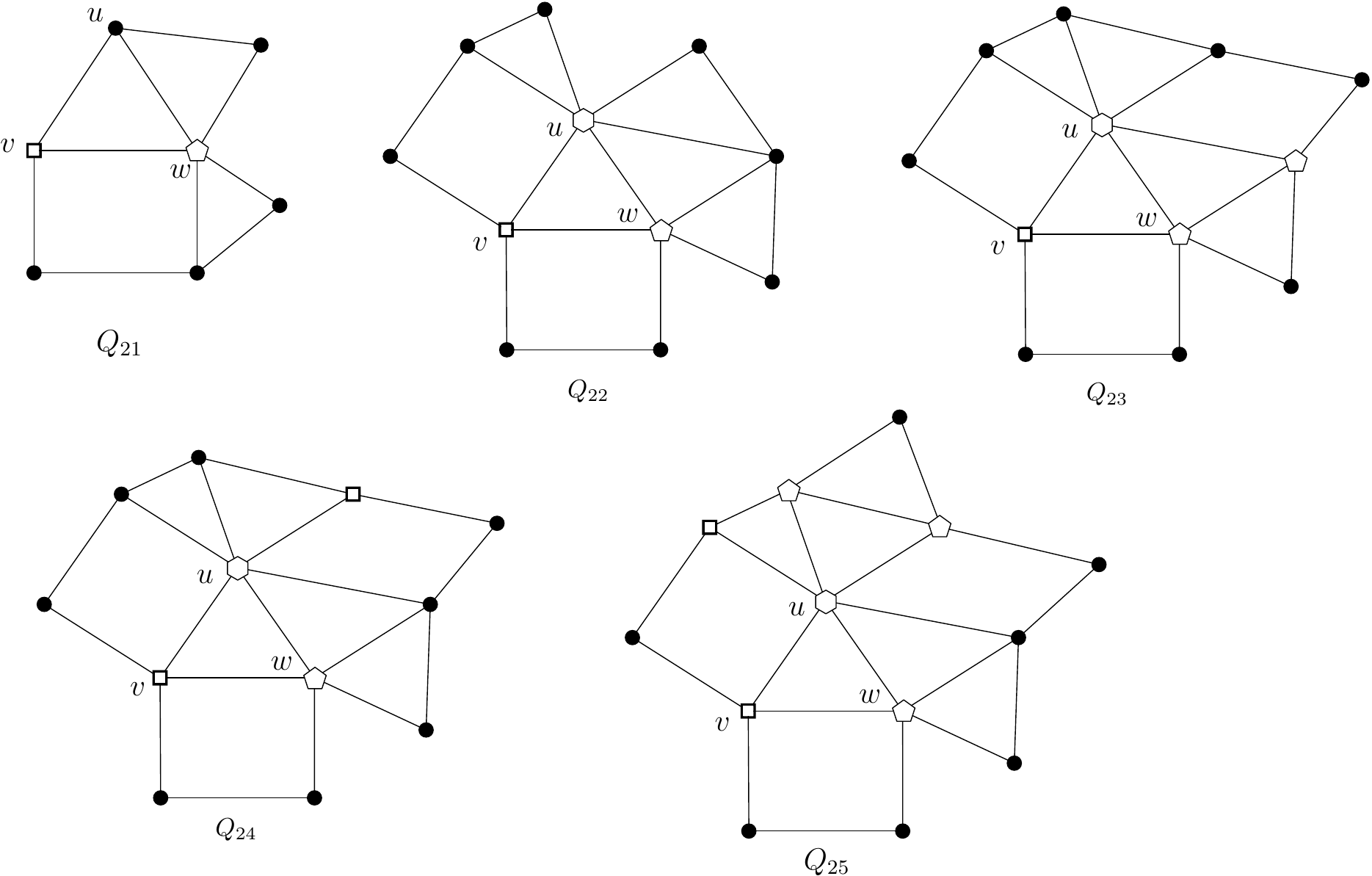}
   \caption{Configurations $Q_{21}$ to $Q_{25}$}
   \label{fig:3.1}
\end{figure}

Configurations used in the paper are shown in Figures \ref{fig:1}--\ref{fig:3.1}.
The vertices shown as squares, pentagons, or hexagons represent the vertices in $U$ and
their values $\delta(u)$ are $4$, $5$, and $6$, respectively. The vertices in $V(C)\setminus U$
are shown as smaller full circles.
The configurations shown in these figures may contains additional
notation that will be used in the proofs later in the paper.

The goal of this section is to prove the following theorem.

\begin{theorem}
\label{thm:unavoidable}
Every plane graph of minimum degree at least four contains one of the
configurations $Q_1,\dots,Q_{25}$ depicted in Figures~\ref{fig:1}--\ref{fig:3.1}.
\end{theorem}

In the proof, we will use the following terminology.
If $v$ is a vertex of degree $k$ in $G$, then we call it a
\DEF{$k$-vertex}, and a vertex of degree at least $k$ (at most
$k$) will also be referred to as a \DEF{$k^{+}$-vertex}
(\DEF{$k^{-}$-vertex}). A neighbor of $v$ whose degree is $k$ is
a \DEF{$k$-neighbor} (similarly \DEF{$k^{+}$-} and
\DEF{$k^{-}$-neighbor}). A face $f$ that has size at least five is
called a \DEF{major face}; if $f$ has size at most 4 it is called
a \DEF{minor face}. A \DEF{$k$-face} is a face of size $k$.
By a \DEF{triangle} we refer to a face of size 3.
An $r$-$s$-$t$ \DEF{triangle} is a triangle whose vertices have
degree $r$, $s$ and $t$, respectively. An $r^{+}$-$s^{+}$-$t^{+}$
\DEF{triangle} is defined similarly.
A triangle is said to be \DEF{bad} if it is a 5-4-4 triangle that is adjacent to at
most two major faces.

\begin{proof}[Proof of Theorem \ref{thm:unavoidable}]
The proof uses the discharging method. Assume, for a
contradiction, that there is a plane graph $G$ that contains none of the configurations shown in
Figures \ref{fig:1}--\ref{fig:3.1}. We shall refer to these configurations as $Q_1,\dots, Q_{25}$.
Let $G$ be a counterexample of minimum order. To each vertex or
face $x$ of $G$, we assign the {\em charge\/} of $c(x)=\deg(x) - 4$. A
well-known consequence of Euler's formula is that the total charge
is always negative, $\sum_{x \in V(G) \cup F(G)} c(x) = -8$. We
are going to apply the following \DEF{discharging rules}:

\begin{itemize}
\item[R1:] A $k$-face ($k \geq 5$) adjacent to $r$ triangles sends
charge of $(k-4)/r$ to each adjacent triangle.

\item[R2:] A 5-vertex $v$ incident to exactly one triangle sends
charge 1 to that triangle. A 5-vertex incident to exactly three
triangles, sends charge $1/3$ to each triangle. A 5-vertex
incident to exactly two triangles sends charge $1/2$ to each
triangle unless (i) at least one of the triangles is a bad
triangle in which case $v$ sends charge of $3/5$ to each bad
triangle and charge of $2/5$ to each non-bad triangle,
or (ii) none of the triangles is bad, one of them is
incident to a 4-vertex and the other is not,
in which case $v$ sends charge $2/3$ to the triangle with the
4-vertex and $1/3$ to the other triangle.

\item[R3:] A 6-vertex $v$ adjacent to a 6-4-4 triangle $T$ sends
charge (i) $4/5$ to $T$ if $T$ is incident to exactly one major
face, (ii) $3/5$ to $T$ if $T$ is incident to exactly two major
faces, and (iii) $2/5$ to $T$ if $T$ is incident to three major
faces.

\item[R4:] A $6^+$-vertex $v$ adjacent to a $6^+$-$5^{+}$-$5^{+}$
triangle $T$ sends charge $1/3$ to $T$ unless $T$ is a $6^+$-5-5
triangle with a $6^+$-5 edge incident to a 4-face and the 5-5 edge
incident to a triangle, in which case $v$ sends charge $7/15$ to $T$.

\item[R5:] A 6-vertex $v$ incident to a 6-5-4 triangle $T=uvw$.
Then $v$ sends charge $1 - x - y$, where $x$ is the total charge
sent to $T$ by the rule R1 and $y$ is the charge sent to $T$ by
the rule R2.

\item[R6:] A 6-vertex $v$ incident to a $6$-$4$-$7^{+}$ triangle
$T$ sends charge $1/3$ to $T$.

\item[R7:] A 6-vertex $v$ incident to a 6-6-4 triangle $T$ sends
charge (i) $1/2$ to $T$ if $T$ is incident to no major faces, (ii)
$2/5$ to $T$ if $T$ is incident to exactly one major face, (iii)
$3/10$ if $T$ is incident to exactly two major faces, (iv) $1/5$
if $T$ is incident to three major faces.

\item[R8:] A $7^{+}$-vertex $v$ incident to a $7^{+}$-4-4 triangle
$T$ sends charge $4/5$ to $T$.

\item[R9:] A $7^{+}$-vertex $v$ incident to a $7^{+}$-$5^{+}$-4
triangle $T$ sends charge $2/3$ to $T$.

\item[R*:] After rules R1--R9 have been applied, each triangle $T$
with positive current charge equally redistributes its excess
charge among those incident 5-5-4 triangles that have negative charge.
\end{itemize}

First, let us state two simple observations that will be used
repeatedly.

\begin{claim}
A $5$-vertex sends charge of at least $1/3$ to every incident triangle.
\end{claim}

\begin{claim}
A major face sends charge of at least $1/5$ to every adjacent triangle.
\end{claim}

For $x \in V(G) \cup F(G)$, let $c^{*}(x)$ be the \DEF{final charge}
obtained after applying rules R1--R9 and R$^*$ to $G$. We will show that
every vertex and face has non-negative final charge. This will
yield a contradiction since the initial total charge of $-8$
must be preserved.

\medskip

\textbf{$4^{+}$-faces:} Since the charge of a $4^+$-face only changes by rule R1,
it is clear that every such face has a nonnegative final charge.

\medskip
\textbf{3-faces:} Let $T = uvw$ be a triangle. Then $c(T) = -1$.
We will show that $c^{*}(T) \geq 0$. We consider a few cases.

\textbf{Case 1: $T$ is a 4-4-4 triangle.} This case is not
possible since $Q_3$ is excluded.

\textbf{Case 2: $T$ is a 5-4-4 triangle.} Let $deg(u)=5 $ and
$deg(v)=deg(w)=4$. We may assume that $u$ is incident to at least
two triangles, for otherwise $c^{*}(T) \geq 0$. Since $Q_1$ and
$Q_5$ are excluded, $u$ is incident to precisely one other
triangle $T'$. If $T$ is not a bad triangle, then all of its
incident faces are major and by R1 and R2, $c^{*}(T) \geq -1 + 1/5
+ 1/5 + 1/5 + 2/5 \geq 0$. Now, suppose that $T$ is a bad
triangle. If $T$ is incident to two major faces, then by R1 and
R2, $c^{*}(T) \geq -1+ 1/5 + 1/5 + 3/5 \geq 0$. Now, assume that
$T$ is incident to at most one major face. Since $Q_1$ and $Q_2$
are excluded, the face incident to the edge $vw$ is the major
face, and the faces incident to $uv$ and $uw$ are both 4-faces.
But now, the exclusion of $Q_9$ implies that $u$ cannot be
incident to any other triangle except $T$, a contradiction.

\textbf{Case 3: $T$ is a $6^{+}$-4-4 triangle.} Let $deg(u) \geq
6$ and $deg(v)=deg(w)=4$. Since $Q_1$ and $Q_2$ are
excluded, $T$ is adjacent to at least one major face. Since a
major face always sends charge at least $1/5$ to an adjacent
triangle, it follows by the rule R3 (if $deg(u)=6$) or R8 (if
$deg(u)\ge7$) that $c^{*}(T) \geq 0$.

\textbf{Case 4: $T$ is a 5-5-4 triangle.}
Let $deg(v)=4$ and $deg(u)=deg(w)=5$. We consider several
subcases.

\textbf{Subcase (a): $T$ is incident to at least two major faces.}
In this case, by rules R1 and R2, $T$ receives total charge of at
least $1/3 + 1/3 + 1/5 + 1/5 > 1$, which implies that $c^{*}(T)\geq 0$.

\textbf{Subcase (b): $T$ is incident to no major faces.} First,
suppose that all faces adjacent to $T$ are 4-faces. Since
$Q_7$ is excluded, each of $u$ and $w$ is incident to at most one
other triangle besides $T$. If $u$ (or $w$) is incident to no
other triangle, then $T$ receives a charge of $1$ from $u$ (or
$w$) and $c^{*}(T) \geq 0$. Therefore, we may assume that each of
$u$ and $w$ is incident to exactly two triangles. But now, the
exclusion of $Q_9$ implies that none of $u$ and $w$ are incident
to a (bad) 5-4-4 triangle. Therefore, each of $u$ and $w$ send
charge of $1/2$ to $T$ by the rule R2. Hence, $c^{*}(T) \geq0$.

The remaining possibility (by exclusion of $Q_1$) is that the face
incident to $uw$ is a triangle and the faces incident to $uv$ and
$vw$ are 4-faces. However, this gives the configuration $Q_6$.

\textbf{Subcase (c): $T$ is incident to exactly one major face.}
We consider several subcases. First, assume that the face incident
to the edge $uw$ is a triangle $T'$. Since $Q_1$ is excluded,
we may assume by symmetry that the face incident to $uv$ is a 4-face
$S$ and the face incident to $vw$ is a major face $R$. If any of
$u$ or $w$ is incident to no other triangles except $T$ and $T'$,
then $T$ receives a total charge of at least $1/2 + 1/3 + 1/5 > 1$
by rules R1 and R2, yielding $c^{*}(T) \geq 0$. Therefore, we may
assume that each of $u$ and $w$ is incident to exactly three
triangles. Since $Q_{21}$ and $Q_{15}$ are excluded, we may assume
that $G$ contains the configuration $P_1$ shown in Figure
\ref{fig:4}. Clearly, $deg(u_3) \geq5$ since $Q_1$ is excluded.
Suppose first that $deg(u_3)=5$. Note that the face $F \neq T'$
incident to the edge $u_3w$ cannot be a triangle since $Q_4$ is
excluded. Now, the exclusion of $Q_{14}$ implies that the face $F$
is a major face. Therefore, $T'$ receives a total charge of at
least $1/3 + 1/3 + 1/3 + 1/5 = 6/5$ after the rules R1--R9 have
been applied. Note that by the exclusion of $Q_{13}$, $T$ is the
only 5-5-4 triangle that is adjacent to $T'$. Hence, $T'$ sends
charge of at least $6/5 - 1= 1/5$ to $T$ by the rule R*.
Therefore, the total charge sent to $T$ is at least $1/3 + 1/3+
1/5 + 1/5 > 1$, resulting in $c^{*}(T)
> 0$. Next, suppose that $deg(u_3) \geq 6$. If $F$ is a major
face, then $T'$ receives a total charge of at least $1/3 + 1/3 +
1/3 + 1/5 = 6/5$ after the rules R1--R9 have been applied. Since
$u_3$ is a $6^{+}$-vertex, $T$ is the only 5-5-4 triangle incident
to $T'$ and therefore by rule R* it receives charge of $1/5$ from
$T'$. As before, $c^{*}(T) \geq 0$. Now, suppose that $F$ is a
4-face. By the rule R4, $u$ sends charge $7/15$ to $T'$.
Hence the total charge received by $T'$ is $1/3+ 1/3+7/15 =
17/15$. This implies that $T'$ sends charge of $17/15 - 1 =2/15$
to $T$ by the rule R*. It follows that $T$ receives total charge
of at least $1/3+1/3+1/5+2/15 = 1$. Hence, $c^{*}(T) \geq 0$.

Next, assume that the face incident to the edge $uw$ is a 4-face.
Since $Q_1$ is excluded, it follows that $T$ is not adjacent to
any triangles. Therefore, we may assume without loss of generality
that the face incident to $uv$ is a 4-face and the face adjacent
to the edge $vw$ is a major face. If none of $u$ and $w$ are
incident to three triangles, then by the rules R1 and R2, $T$
receives a charge of at least $2/5 + 2/5 + 1/5 = 1$, and we are
done. Therefore, we may suppose that either $u$ or $w$ is incident
to three triangles. First, suppose that $u$ is incident to three
triangles. By the exclusion of $Q_{16}$ and $Q_{17}$, $w$ is
incident to at most two triangles and no triangle incident to $w$
is a 5-4-4 triangle. By the rule R2, $w$ sends charge of at least
$1/2$ to $T$. Hence, $T$ receives total charge of at least $1/2+
1/3+ 1/5 > 1$, which implies that $c^{*}(T) \geq 0$. Now, suppose
that $w$ is incident to three triangles. By the exclusions of
$Q_{16}$ and $Q_{18}$, $u$ is incident to at most two triangles
none of which is a 5-4-4 triangle. Therefore, $u$ sends charge of
at least $1/2$ to $T$, and as before we obtain that $c^{*}(T)
\geq0$.

Finally, assume that the face incident to the edge $uw$ is a major
face. Since $Q_1$ is excluded it follows that the faces incident
to the edges $uv$ and $vw$ are both 4-faces. If none of $u$ and
$w$ are incident to three triangles, then by the rules R1 and R2,
$T$ receives a charge of at least $2/5 + 2/5 + 1/5 = 1$, and we
are done. Therefore, we may assume by symmetry
that $u$ is incident to three triangles. By the exclusion of
$Q_{19}$ and $Q_{20}$, $w$ is incident to at most two triangles,
none of which is a 5-4-4 triangle. By R2, $w$ sends charge of at
least $1/2$ to $T$. Hence, $T$ receives total charge of at least
$1/2+ 1/3+ 1/5 > 1$, which implies that $c^{*}(T) \geq 0$.

\begin{figure}[htb]
   \centering
   \includegraphics[width=13cm]{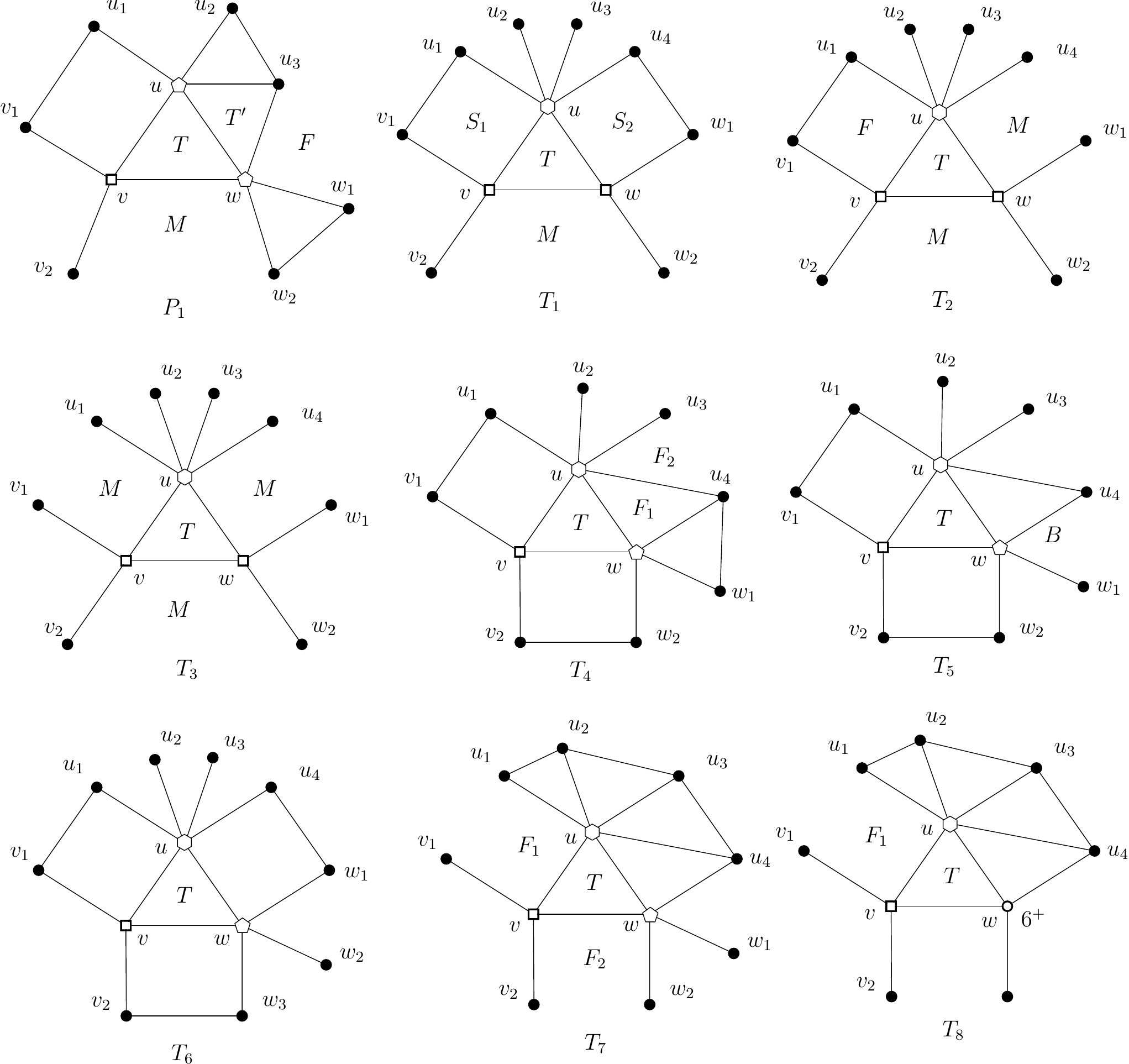}
   \caption{Discharging analysis}
   \label{fig:4}
\end{figure}

\textbf{Case 5: $T$ is a 6-5-4 triangle.}
By the rules R1, R2 and R5, $T$ receives a total charge of 1 after
the discharging rules have been applied. Hence, $c^{*}(T) \geq 0$.

\textbf{Case 6: $T$ is a 6-6-4 triangle.}
Clearly, only rules R1 and R7 apply to $T$. If $T$ is adjacent to
no major faces, each 6-vertex contained in $T$ sends charge of
$1/2$ to $T$, which yields $c^{*}(T) \geq 0$. If $T$ is adjacent
to one major face, each 6-vertex contained in $T$ sends charge of
$2/5$ to $T$, and $T$ receives a charge of at least $2/5 + 2/5 +
1/5=1$. If $T$ is adjacent to two major faces, each 6-vertex
contained in $T$ sends charge of $3/10$ to $T$, and $T$ receives
charge of at least $3/10+3/10+ 1/5 + 1/5=1$, and we are done. If
$T$ is adjacent to three major faces, each 6-vertex contained in
$T$ sends charge of $1/5$ to $T$, and $T$ receives charge of at
least $1/5+1/5+ 1/5 + 1/5 + 1/5=1$. Thus, in all cases $c^{*}(T)
\geq 0$.

\textbf{Case 7: $T$ is a $7^{+}$-$5^{+}$-4 triangle.}
By rule R9, the $7^{+}$-vertex incident to $T$ sends charge of
$2/3$ to $T$. By rules R2, R6 and R9, the $5^{+}$-vertex sends
charge of at least $1/3$ to $T$. Thus, $T$ receives a total charge
of at least $1$, and $c^{*}(T) \geq 0$.

\textbf{Case 8: $T$ is a $5^{+}$-$5^{+}$-$5^{+}$ triangle.}
Since the rule R7 is clearly not applied to $T$, each
$5^{+}$-vertex incident to $T$ sends charge of at least $1/3$ to
$T$. Hence, $c^{*}(T) \geq 0$.

\textbf{4-vertices:} If $v$ is a 4-vertex, then $c(v)=0$. Since
$v$ neither receives nor gives any charge, we have $c^{*}(v)=0$.

\textbf{5-vertices:} Let $v$ be a 5-vertex. If $v$ is incident to
exactly one triangle, then clearly by rule R2, $c^{*}(v)=0$. If
$v$ is incident to exactly three triangles, again by rule R2,
$c^{*}(v)=0$. Note that $v$ is incident to at most three triangles
since $Q_4$ is excluded. Hence, it remains to consider the case
that $v$ is incident to exactly two triangles $T_1$ and $T_2$. To
show that $c^{*}(v) \geq 0$, by R2 it is sufficient to show that
not both of $T_1$ and $T_2$ are bad. This follows from the
exclusion of $Q_1$, $Q_2$, $Q_8$ and $Q_9$.

\textbf{6-vertices:} A 6-vertex $u$ has initial charge of 2.
We break the analysis into several cases,
depending on what type of triangles $u$ is incident to. Often we
will reduce the analysis to previously considered cases.

\textbf{Case 1: $u$ is incident to a 6-4-4 triangle.} Suppose
$T=uvw$ is a 6-4-4 triangle, with $deg(v)=deg(w)=4$. Since $Q_1$
and $Q_2$ are excluded, none of adjacent faces to $T$ are
triangles, and the face incident to $vw$ is a major face. The
possible cases are outlined in Figure \ref{fig:4} as
configurations $T_1$--$T_3$. When such a configuration $T_i$
occurs at $u$, we say that $u$ is a \DEF{type $T_i$ vertex}. Note
that in the figure, a face marked by $M$ stands for a major face.

\textbf{Subcase (a): $u$ is a type $T_1$ vertex.} Note that $u$
sends charge of 4/5 to $T$ by rule R3. We consider a few cases
depending on the number of triangles incident to $u$. Since
$Q_{11}$ is excluded, $u$ cannot be incident to four triangles.
Since a 6-vertex never sends charge more than $4/5$ to a triangle,
$u$ sends total charge of at most $8/5 < 2$, yielding a positive
final charge for $u$ when $u$ is incident to at most two
triangles. Therefore, we may assume that $u$ is incident to
exactly three triangles. Since $Q_{11}$ is excluded, we may assume
without loss of generality that $u_2u_3, u_3u_4$ are edges, but
$u_1u_2$ is not. Clearly, by exclusion of $Q_1$, $u_3$ is a
$5^{+}$-vertex. Therefore, $u$ sends charge to the triangles
$uu_2u_3$ and $uu_3u_4$ by rules R4--R7. The most charge sent out
by rule R5 is $2/3$ and most sent out by other rules is $1/2$.
Since $4/5+2/3+1/2 < 2$ and $c(u)=2$, it suffices to check what
happens if rule R5 is applied twice. In such a case we have
$deg(u_2)=deg(u_4)=4$ and $deg(u_3)=5$. Since $Q_{10}$ is
excluded, the face $F$ bounded by $u_1u$ and $uu_2$ is a major
face. But in this case, by rule R5, $u$ sends charge of at most
$1-1/3-1/5=7/15$ to the triangle $uu_2u_3$. Thus, the final charge
sent by $u$ is again at most $7/15 + 2/3 + 4/5 < 2$.

\textbf{Subcase (b): $u$ is a type $T_2$ vertex.} First, assume
that $u$ is incident to four triangles. This implies that
$u_1u_2$, $u_2u_3$ and $u_3u_4$ are all edges. Suppose that $u_1$
is a $5^{+}$-vertex. By exclusion of $Q_1$, $u_2$ and $u_3$ are
$5^{+}$-vertices. By rule R4, $u$ sends charge at most $7/15$ to
the triangle $uu_1u_2$, and $1/3$ to the triangle $u_2uu_3$. Since
the triangle $uu_3u_4$ is incident to at least one major face, by
rules R4-R7, $u$ sends charge of at most $7/15$ to the triangle
$uu_3u_4$. By rule R3, $u$ sends charge $3/5$ to $T$. Therefore,
the total charge sent by $u$ is at most $7/15 + 1/3 + 7/15 + 3/5 =
28/15 < 2$, yielding $c^{*}(u) > 0$. On the other hand, if $u_1$
is a 4-vertex, exclusion of $Q_{12}$ implies that $deg(u_2) \geq
6$. Again, we see that $u$ sends charge of at most $1/2 + 1/3 +
7/15 + 3/5 < 2$.

Now, assume that $u$ is incident to at most three triangles. By
subcase (a), we may assume that $u$ is not a type $T_1$ vertex
(for another triangle incident with $u$).
Hence, we may assume $u$ never sends charge $4/5$ to a triangle.
By rules R3--R7, it is clear to see that $u$ sends charge at most
$2/3$ to each triangle. Since there are at most three triangles,
$u$ sends total charge of at most 2, resulting in $c^{*}(u) \geq0$.

\textbf{Subcase (c): $u$ is a type $T_3$ vertex.} Again, $u$ is
incident to at most four triangles. First, assume that $u$ is
incident to exactly four. This implies that $u_1u_2$, $u_2u_3$ and
$u_3u_4$ are all edges. By the exclusion of $Q_1$, $u_2$ and $u_3$
are both $5^{+}$-vertices. This implies that $u$ sends charge
$1/3$ to the triangle $uu_2u_3$.

Since the triangle $uu_1u_2$ is incident to at least one major
face, by rules R4-R7, $u$ sends charge of at most $7/15$ to the
triangle $uu_1u_2$. Similarly, $u$ sends at most $7/15$ to
the triangle $uu_3u_4$. By rule R3, $u$ sends charge at most $2/5$
to the triangle $T$. Therefore, $u$ sends total charge of at most
$1/3 + 7/15 + 7/15 + 2/5 < 2$, yielding $c^{*}(u) \geq 0$.

The case that $u$ is incident to at most three triangles is
handled by an identical argument as in the previous subcase.

\textbf{Case 2: $u$ is incident to a 6-5-4 triangle.} Suppose
$T=uvw$ is a 6-5-4 triangle, with $deg(v)=4$ and $deg(w)=5$. Since
$Q_1$ is excluded, the faces different from $T$ that are incident to the
edges $uv$ and $vw$ are not triangles. We consider the subcases shown
in Figure \ref{fig:4} as $T_4$--$T_7$; in $T_7$, at least one of $F_1,F_2$
is a major face, and it will be argued later why we may assume that $u$ is
incident with 5 triangles.
We may assume that $u$ is not incident to any 6-4-4 triangle, and
consequently, never sends more than $2/3$ charge to any triangle.
In particular, we may assume that $u$ is incident with four or
five triangles.

\textbf{Subcase (a): $u$ is a type $T_4$ vertex.}
In this case, $u$ is incident to precisely four triangles since the configuration
$Q_{22}$ is excluded. Additionally, $uu_2u_3$ is a triangle.
First, assume that $u_3u_4$ is an edge (and consequently, $u_1u_2$ is not an
edge).  Since $Q_1$ is excluded, we have that both $u_3$ and $u_4$
are $5^{+}$-vertices. Therefore, by rule R4, $u$ sends charge
$1/3$ to each of the triangles $uwu_4$ and $uu_3u_4$. Since $u$
sends charge of at most $2/3$ to each of the other two incident
triangles, we have $c^{*}(u) \geq 0$.

Second, assume that $u_1u_2$ is an edge (consequently, $u_3u_4$ is
not an edge). We divide the analysis into two cases, depending on
whether the face $F_2$ bounded by the edges $uu_3$ and $uu_4$ is a
4-face or a major face. If $F_2$ is a 4-face, then the exclusion
of $Q_{23}$ implies that $u_4$ is a $6^{+}$-vertex, and the
exclusion of $Q_{24}$ implies that $u_3$ is a $5^{+}$-vertex.
Therefore, $u$ sends charge of $1/3$ to the triangle $F_1$ by rule
R4. If $u$ sends charge of at most $1/3$ to another incident
triangle, then $c^*(u)\ge0$. Thus, we may assume that $u$
sends charge $7/15$ to the triangle $uu_2u_3$ by the last subcase
in rule R4. Therefore, $\deg(u_2)=\deg(u_3)=5$. By excluding
$Q_{25}$, we see that $u_1$ is a $5^{+}$-vertex. By rule R4, $u$
sends total charge of at most $14/15 < 1$ to the triangles
$uu_1u_2$ and $uu_2u_3$, and we are done.

Now, suppose that $F_2$ is a major face. Then $u$ sends charge
$1/3$ to the triangle $uu_4w$ by rule R4. Since $u$ sends charge
of $2/3$ to $T$, it suffices to show that $u$ sends total charge
of at most $1$ to the triangles $uu_1u_2$ and $uu_2u_3$. If $u_2$
is a $6^{+}$-vertex, then $u$ sends charge of at most $2/5$ to the
triangle $uu_2u_3$ and charge of at most $1/2$ to the triangle
$uu_1u_2$ by one of the rules R4, R6 or R7. Otherwise, $u_2$ is a
5-vertex. Then $u$ sends at most $7/15$ to the triangle
$uu_2u_3$. If $u_1$ is a $5^{+}$-vertex, $u$ sends charge at most
$7/15$ to the triangle $uu_1u_2$ by the rule R4, and we are done.
Therefore, we may assume that $u_1$ is a 4-vertex. We may assume
that the second face incident to the edge $u_1u_2$ is a 4-face, for
otherwise $u$ sends charge at most $7/15$ to the triangle
$uu_1u_2$ by rule R5, and we are done. Furthermore, if $u_2$ is
not incident to any other triangles except $uu_1u_2$ and
$uu_2u_3$, then by rule R5, $u$ sends charge of at most $1/2$ to
the triangle $uu_1u_2$ and again, we are done. Therefore, we may
assume that $u_2$ is incident to a third triangle $T'$. Since the
configuration $Q_{21}$ is forbidden, $T'$ contains the edge
$u_2u_3$. But now, the exclusion of $Q_1$ implies that $u_3$ is a
$5^{+}$-vertex, and in fact, $u$ sends charge of $1/3$ to the
triangle $uu_2u_3$. Since $u$ sends charge of at most $2/3$ to the
triangle $uu_1u_2$, we are done.

\textbf{Subcase (b): $u$ is a type $T_5$ vertex.} Note that $B$ in
the figure denotes a ``big'' face -- a face of size at least 4.
Also, note that the face determined by the edges $ww_1$ and
$ww_2$ cannot be a triangle since $Q_{21}$ is excluded. Therefore,
in this subcase $w$ is incident to two triangles. Therefore, $w$
sends charge $2/3$ to $T$ (by rule R2(ii)) if $deg(u_4)\ne4$.
Thus, if $deg(u_4)\ne4$, $u$ sends $1/3$ to $T$.

First, assume that $u$ is incident to five triangles, i.e., the
edges $u_1u_2, u_2u_3, u_3u_4$ are all present. The
exclusion of $Q_1$ implies that $u_2, u_3, u_4$ are all
$5^{+}$-vertices. Now, by rule R5, $u$ sends charge of $1/3$ to
$T$. By rule R4, $u$ sends charge of $1/3$ to each of the
triangles $uu_4w$, $uu_3u_4$ and $uu_2u_3$. By rules R4--R7, $u$
sends charge of at most $2/3$ to the triangle $uu_1u_2$.
Therefore, the total charge sent by $u$ is at most $1/3 + 1/3 +
1/3+ 1/3 + 2/3=2$, which implies that $c^{*}(u) \geq 0$.

Next, assume that $u$ is incident to exactly four triangles. There
are three possibilities.

First, assume that $u_2u_3$ and $u_3u_4$ are edges (and $u_1u_2$
is not). Then the exclusion of $Q_1$ implies that $u_3$ and $u_4$
are both $5^{+}$-vertices. As before, $u$ sends charge of $1/3$ to
$T$, charge of $1/3$ to $uu_4w$ and $1/3$ to $uu_3u_4$. Since $u$
sends charge of at most $2/3$ to the triangle $uu_2u_3$, we have
$c^*(u)\ge0$.

Second, assume that $u_1u_2$ and $u_3u_4$ are edges (and $u_2u_3$
is not). As before, $u_4$ is a $5^{+}$-vertex. Then $u$ sends
charge $1/3$ to each of the triangles $T$ and $uwu_4$, and charge
of at most $2/3$ to each of the triangles $uu_1u_2$ and $uu_3u_4$.
Thus, $u$ sends total charge of at most $1/3+ 1/3 + 2/3+ 2/3 = 2$.

Finally, assume that $u_1u_2$ and $u_2u_3$ are edges, and $u_3u_4$
is a non-edge. If $u_4$ is a 4-vertex, then $u$ sends charge of
$1/2$ to $T$ and charge of at most $1/2$ to the triangle $uu_4w$.
If $u_4$ is a $5^{+}$-vertex, then $u$ sends charge of $1/3$ to
$T$ and charge of at most $7/15$ to $uu_4w$. Thus, $u$ sends
charge of at most 1 to the triangles $T$ and $uu_4w$. Therefore,
it is sufficient to show that $u$ sends a total charge of at most
1 to the triangles $uu_1u_2$ and $uu_2u_3$.

Clearly, $u_2$ is a $5^{+}$-vertex since $Q_1$ is excluded. If
$u_2$ is a $6^{+}$-vertex, then $u$ sends charge at most $1/2$ to
each of the triangles $uu_1u_2$ and $uu_2u_3$ by the rules R4, R6
and R7. Therefore, we may assume that $u_2$ is a 5-vertex. It is
sufficient to show that $u$ sends charge of at most $1/2$ to each
of the triangles $uu_1u_2$ and $uu_2u_3$. If $u_1$ is a
$5^{+}$-vertex then by rule R4, $u$ sends charge of at most $7/15$
to $uu_1u_2$. If $u_1$ is a 4-vertex, then by the exclusion of
$Q_{21}$ (and $Q_1$), the triangle $uu_1u_2$ is either incident to
a major face, or $u_2$ is incident to only two triangles or $u$ is
a type $T_4$ vertex (for the triangle $uu_1u_2$).
Since we may assume that $u$ is not a type
$T_4$ vertex, it follows that $u$ sends charge of at most $1/2$ to
the triangle $uu_1u_2$. A similar argument applied to $u_3$ shows
that $u$ sends charge of at most $1/2$ to the triangle $uu_2u_3$.

\textbf{Subcase (c): $u$ is a type $T_6$ vertex.} Since we may
assume that $u$ is not a type $T_1$ vertex, by rules R3--R7, $u$
never sends charge of more than $2/3$ to an incident triangle.
Therefore, if $u$ is incident to at most three triangles, we have
that $c^{*}(u) \geq 0$. Since $u$ is a type $T_6$ vertex, it is
incident to four triangles. Therefore, the only possibility we
have is when $u_1u_2, u_2u_3, u_3u_4$ are all edges. By exclusion
of $Q_1$, we have that $u_2$ and $u_3$ are both $5^{+}$-vertices.
It follows that $u$ sends charge of $1/3$ to the triangle
$uu_2u_3$ by rule R4. Since $u$ sends charge of at most $2/3$ to
the triangle $T$, it is sufficient to show that $u$ sends charge
of at most $1/2$ to each of the triangles $uu_1u_2$ and $uu_3u_4$.
Consider the triangle $uu_1u_2$. If $u_2$ is a $6^{+}$-vertex,
then by rules R4, R6 and R7, $u$ sends charge of at most $1/2$ to
the triangle $uu_1u_2$. Therefore, we may assume that $u_2$ is a
5-vertex. We may assume that $u_1$ is a 4-vertex, for otherwise by
rule R4, $u$ sends charge of at most $7/15$ to $uu_1u_2$, and we
are done. This implies (by the exclusion of $Q_1$) that the second
face incident to the edge $u_1u_2$ is a $4^{+}$-face. If it is a
$4$-face, then $u$ is actually a type $T_4$ or $T_5$ vertex, and
we are done by the previous analysis. Therefore, we may assume
that the edge $u_1u_2$ is incident to a major face. But then, by
the rule R5, $u$ sends charge of at most $7/15$ to the triangle
$uu_1u_2$. An identical argument shows that $u$ sends charge of at
most $1/2$ to the triangle $uu_3u_4$. Thus, in all cases $c^{*}(u)
\geq 0$.

\textbf{Subcase (d): $T$ is incident to at least one major face.}
Since we may assume that $u$ is not a type $T_i$ vertex, for any
$1 \leq i \leq 6$, it follows that $u$ is not incident to a 6-4-4
triangle, and sends a charge of at most $7/15$ to any incident
6-5-4 triangle. Therefore, by rules R4--R7, $u$ sends a charge of
at most $1/2$ to any incident triangle. It follows that if $u$ is
incident to at most four triangles, then $c^{*}(u) \geq 0$. Since
the edge $uv$ is incident to a $4^{+}$-face by the exclusion of
$Q_1$, we may assume that we have the configuration $T_7$ of
Figure \ref{fig:4}, where at least one of the faces $F_1$ and
$F_2$ is a major face. It follows that the vertices $u_2, u_3,
u_4$ are all $5^{+}$ vertices, and by rule R4, $u$ sends charge of
$1/3$ to each of the triangles $uu_2u_3$, $uu_3u_4$ and $uu_4w$.
Since either $F_1$ or $F_2$ is a major face, it follows that $u$
sends charge of at most $7/15$ to $T$ by the rule R5. Since $u$
sends charge of at most $1/2$ to any triangle, we have that the
total charge sent by $u$ is at most $1/3+ 1/3 + 1/3 + 7/15 + 1/2 <
2$, yielding $c^{*}(u)>0$.

\textbf{Case 3: $u$ is incident to a $6$-$6^{+}$-$4$ triangle.}
Suppose $T=uvw$ is a 6-$6^{+}$-4 triangle, with $deg(v)=4$. By
cases 1 and 2, we may assume that $u$ is neither incident to a
6-4-4 triangle nor to a 6-5-4 triangle. Therefore, rules R3 and R5
never apply to $u$, and by rules R4, R6 and R7, $u$ sends charge
of at most $1/2$ to each incident triangle. Thus, if $u$ is
incident to at most four triangles, $c^{*}(u) \geq 0$. Therefore,
we may assume that $u$ is incident to at least five triangles.
Since the face $F_1$ incident to the edge $uv$ cannot be a
triangle (by the exclusion of $Q_1$), it follows that the only
possibility left to consider is the configuration $T_8$ in Figure~
\ref{fig:4}. Now, by exclusion of $Q_1$, $u_2, u_3, u_4$ are all
$5^{+}$ vertices, and hence, by rule R4, $u$ sends charge of $1/3$
to each of the triangles $uu_2u_3$, $uu_3u_4$ and $uu_4w$. Since
$u$ never sends charge of more than $1/2$ to an incident triangle,
we get that $u$ sends a total charge of at most $1/3+1/3+1/3+ 1/2+
1/2 = 2$, as required.

\textbf{Case 4: $u$ is incident to a $6$-$5^{+}$-$5^{+}$ triangle.}
Suppose $T=uvw$ is a 6-$5^{+}$-$5^{+}$ triangle. By cases 1--3,
we may assume that $u$ is not incident to any triangle that
contains a 4-vertex. Therefore, only rule R4 applies to $u$, and
consequently $u$ sends charge of either $1/3$ or $7/15$ to any
incident triangle. Therefore, if $u$ is incident to at most four
triangles, $c^{*}(u) \geq 0$. If $u$ is incident to six triangles,
then rule R4 implies that $u$ sends charge $1/3$ to each triangle,
yielding $c^{*}(u) = 0$. Therefore, we may assume that $u$ is
incident to exactly five triangles. But then it is clear that
there are at most two triangles incident to $u$ to which it sends
charge of $7/15$. Thus, in this case as well, $u$ sends charge of
at most $7/15 + 7/15+ 1/3 + 1/3 + 1/3 < 2$, as required.

\textbf{$7^{+}$-vertices:} Let $u$ be a $7^{+}$-vertex. First,
assume that $deg(u) = d \geq 8$. Note that by rules R4, and
R8--R9, $u$ sends charge of at most $4/5$ to any incident
triangle, and charge of $2/3$ if it is not incident to an
$8^{+}$-4-4 triangle. Therefore, if $u$ is incident to at most
$d-3$ triangles, then it sends total charge of at most
$\tfrac{4}{5}(d-3) \leq d-4$, since $d\geq 8$. Therefore, in this
case $c^{*}(u) \geq 0$. Hence, we may assume that $u$ is incident
to at least $d-2$ triangles. If $u$ is incident to at least $d-1$
triangles, then the exclusion of $Q_1$ implies that at least $d-3$
of these triangles are $8^{+}$-$5^{+}$-$5^{+}$ triangles and that
none of the triangles is an $8^{+}$-4-4 triangle. Note that none
of the $d-3$ mentioned $8^{+}$-$5^{+}$-$5^{+}$ triangles have a
4-face incident to an $8^{+}$-$5^{+}$ edge. Hence, $u$ sends
charge of at most $\tfrac{1}{3}(d-3) + 3 \left( \tfrac{2}{3}
\right) \leq d-4$ for $d \geq 8$. Hence, in this case as well
$c^{*}(u) \geq 0$. Thus, we may assume that $u$ is incident to
exactly $d-2$ triangles. Now, the exclusion of $Q_1$ implies that
$u$ is incident to at most one $8^{+}$-4-4 triangle and at least
one $8^{+}$-$5^{+}$-$5^{+}$ triangle. Therefore, $u$ sends charge
of at most $4/5 + 7/15 + \tfrac{2}{3}(d-4) \leq d-4$, and
$c^{*}(u) \geq 0$.

Now, suppose that $d=7$. Then $u$ has charge +3. Since $u$
sends charge of at most $4/5$ to any incident triangle, we may
assume that $u$ is incident to at least four triangles. If $u$ is
incident to seven triangles, then by exclusion of $Q_1$ they are
all $7$-$5^{+}$-$5^{+}$ triangles, and hence $u$ sends total
charge of $7/3 < 3$. If $u$ is incident to six triangles, then the
exclusion of $Q_1$ implies that $u$ is incident to at least four
$7$-$5^{+}$-$5^{+}$ triangles and to no 7-4-4 triangle. Clearly,
$u$ sends charge of at most $1/3$ to each $7$-$5^{+}$-$5^{+}$
triangle. Therefore, in this case $u$ sends total charge of at
most $4 \left( \tfrac{1}{3} \right) + 2 \left( \tfrac{2}{3}
\right) < 3$. Now, suppose that $u$ is incident to five triangles.
We consider two cases. First, suppose that $u$ is incident to a
7-4-4 triangle. Note that by exclusion of $Q_1$, $u$ is incident
to at most one such triangle. Also, by exclusion of $Q_1$, we have
that $u$ is incident to at least two 7-$5^{+}$-$5^{+}$ triangles.
Clearly, $u$ sends charge of $1/3$ to each of these
7-$5^{+}$-$5^{+}$ triangles. Therefore, $u$ sends total charge of
at most $4/5 + 1/3 + 1/3 + 2/3 + 2/3 < 3$. Secondly, suppose that
$u$ is incident to no 7-4-4 triangle. Then $u$ sends charge of at
most $2/3$ to any incident triangle. Since $u$ is incident to five
triangles, the exclusion of $Q_1$ implies that $u$ is incident to
at least one $7$-$5^{+}$-$5^{+}$ triangle to which it only sends
charge of $1/3$. Therefore, $u$ sends total charge of at most
$\tfrac{1}{3}+ 4 \left( \tfrac{2}{3} \right) = 3$, which implies
that $c^{*}(u) \geq 0$. Lastly, suppose that $u$ is incident to
four triangles. This implies by the exclusion of $Q_1$ that $u$ is
incident to at most two 7-4-4 triangles. Hence, $u$ sends total
charge of at most $4/5 + 4/5 + 2/3 + 2/3 = 44/15 < 3$. Thus, in
all cases, $c^{*}(u) \geq 0$.
\end{proof}

\section{Reducibility}

Let us first introduce some notation that will be used in the rest of the paper. If either $uv$ or $vu$ is an arc in a digraph $D$,
we say that $uv$ is an \DEF{edge} of $D$. We will consider a planar digraph $D$, its underlying graph $G$, and a
2-list-assignment $L$, where $L(v) \subseteq \C$ and $|L(v)|=2$ for every $v\in V(D)$. Given a non-proper $L$-coloring $\phi$ of $D$,
a \DEF{color-$i$ cycle} is a directed cycle in $D$ whose every
vertex is colored with color $i$, for $i \in \C$.
When we speak of vertex \DEF{degrees}, we always mean
degrees in $G$. For the digraph $D$, the \DEF{out-degree} and the
\DEF{in-degree} of a vertex $v$ are denoted by $d^{+}(v)$ and
$d^{-}(v)$, respectively.

If $D$ is a digraph drawn in the plane and $C$ is a configuration (which is an undirected graph), we say that $D$ \DEF{contains} the configuration $C$
if the underlying undirected graph $G$ of $D$ contains $C$.
A configuration $C$ is called \DEF{reducible} if it cannot occur
in a minimum counterexample to Theorem \ref{thm:main choosable}. Showing
that every planar digraph $D$ of minimum degree at least 4 and with digirth at least five contains a reducible configuration
will imply that every such digraph is
2-choosable.

Throughout this section, we assume that $D$ is a planar digraph with digirth at least five that is a counterexample to the theorem
with a 2-list-assignment $L$ such that every proper subdigraph of
$D$ is $L$-colorable. In most statements, we will consider a
special vertex $v \in V(D)$, and we shall assume that $L(v) = \{1,2\}$.
The following lemma shows that the minimum degree of $D$ is at least four and that each vertex has in-degree and out-degree at least two.

\begin{lemma}
\label{lem:min in-out}
Let $v \in V(D)$. Then in every $L$-coloring of $D-v$,
each color in $L(v)$ appears at least once among the out-neighbors and
at least once among the in-neighbors of~$v$.
Consequently, every $v \in V(D)$ has $d^{+}(v) \geq 2$ and
$d^{-}(v)\geq 2$; therefore, $D$ contains no $3^{-}$-vertices
and every $4$-vertex has $d^{+}(v)=d^{-}(v)=2$.
\end{lemma}

\begin{proof}
Suppose that a color $c\in L(v)$ does not appear among the outneighbors of $v$ in an $L$-coloring of $D-v$.
Then coloring $v$ with $c$ gives an $L$-coloring of $D$ since a color-$c$ cycle would have to use an outneighbor of $v$.
The same contradiction is obtained if a color in $L(v)$ does not occur among the in-neighbors, and this completes the proof.
\end{proof}

Having an $L$-coloring $\phi$ of a subdigraph $D-u$ ($u\in V(D)$), we may consider coloring $u$ with a color $i\in L(u)$.
Since $D$ is not $L$-colorable, this creates a color-$i$ cycle, which we denote by $C_i=C_i(u)$.
Such cycles will always be taken with respect to a partial coloring $\phi$ that will be clear from the context.
If $L(u)=\{a,b\}$, then $C_a(u)$ and $C_b(u)$ are disjoint apart from their common vertex $u$.
Since $D$ is drawn in the plane, these cycles cannot cross each other at $u$, and we say that they \DEF{touch}.

\begin{lemma} \label{lem:directed triangle}
Let $v$ be a vertex incident to a triangle $T=vwu$, let $\phi$
be an $L$-coloring of $D-v$, and let $i\in L(v)$.
Then $C_i(v)$ cannot contain both edges $vu$ and $vw$.
\end{lemma}

\begin{proof}
Since $C_i(v)$ is directed, we may assume that $uv, vw \in E(D)$. Since
$D$ has digirth greater than three, this implies that $uw \in E(D)$.
But then we have a color-$i$ cycle in $D-v$ consisting of the path $C_i(v)-v$ and the arc $uw$,
a contradiction.
\end{proof}

\begin{lemma} \label{lem:directed square}
Let $v$ be a vertex incident to a 4-cycle $T=vwux$, let $\phi$
be an $L$-coloring of $D-v$ and let $i\in L(v)$.
Then $C_i(v)$ cannot contain all three edges $ux$, $xv$ and $vw$.
\end{lemma}

\begin{proof}
Suppose that $C_i(v)$ contains the edges $ux$, $xv$ and
$vw$. Since $C_i(v)$ is directed, we may assume that $ux, xv, vw \in
E(D)$. But this implies that $uw \in E(D)$, and we have a color-$i$
cycle through the arc $uw$ in $D-v$, a contradiction.
\end{proof}

\begin{figure}[htb]
   \centering
   \includegraphics[width=14cm]{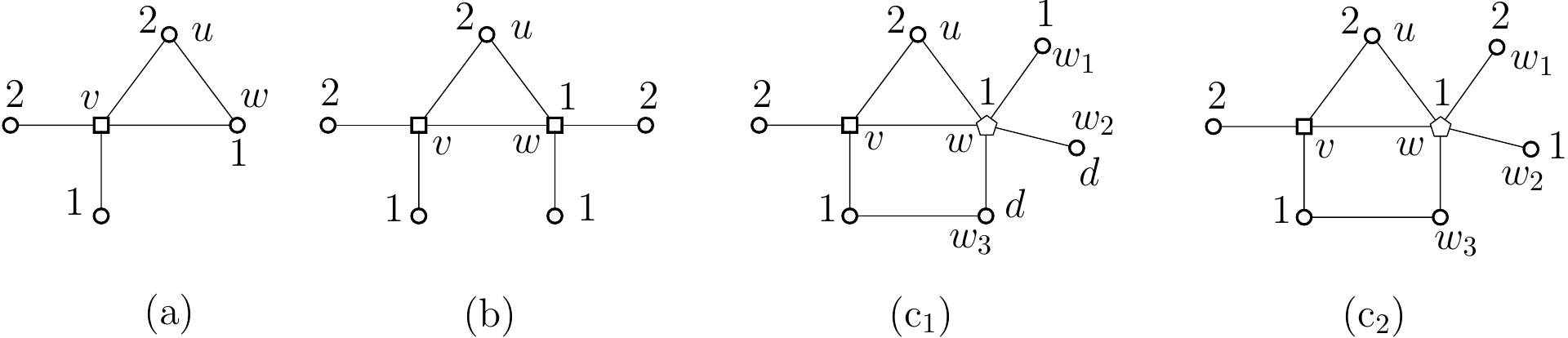}
   \caption{Colors around a 4-vertex contained in a triangle}
   \label{fig:colors}
\end{figure}

The next lemma shows some restrictions on the colors around a 4-vertex that is contained in a triangle.
Recall our assumption that $L(v)=\{1,2\}$.

\begin{lemma}
\label{lem:triangle}
Let $T=vuw$ be a triangle in $D$ and $deg(v)=4$.
Let $\phi$ be an $L$-coloring of $D-v$ such that $\phi(w)=1$.

{\rm (a)} The colors of the neighbors of $v$ and the cycles $C_1(v)$ and $C_2(v)$ are as shown in Figure \ref{fig:colors}(a).

{\rm (b)} If $deg(w)=4$, then $L(w)=L(v)=\{1,2\}$ and the colors of the neighbors of $v$ and $w$ are as shown in Figure \ref{fig:colors}(b).

{\rm (c)} If $deg(w)=5$, the other face containing the edge $vw$ is a 4-face, and the clockwise neighbors of $w$ are $v,u,w_1,w_2,w_3$,
then either (i) $w_1\in V(C_1(v))$ and the colors of the neighbors of $v$ and $w$ are as shown in
Figure \ref{fig:colors}(c$_1$), where $d$ is the color in $L(w)\setminus \{1\}$), or 
(ii) $L(w)=L(v)=\{1,2\}$, and the colors of the neighbors of $v$ and $w$ are as shown in
Figure \ref{fig:colors}(c$_2$).
\end{lemma}

\begin{proof}
(a) By assumption, $w\in V(C_1(v))$. By Lemma \ref{lem:directed triangle}, $u\notin V(C_1(v))$. Since $C_1(v)$ and $C_2(v)$ touch at $v$,
they must be as claimed.

(b) By uncoloring $w$ and coloring $v$ with color 1, we obtain an $L$-coloring $\phi'$ of $G-w$. The claim follows by applying
part (a) to $D-w$ and $\phi'$.

(c) By (a), colors around $v$ are as claimed. By Lemma \ref{lem:directed square}, the cycle $C_1(v)$ does not contain $w_3$.
We are done if it contains $w_1$. Thus, we may assume that $C_1(v)$ contains $w_2$. Let us consider the coloring $\phi'$ of $D-w$
as used in the proof of part (b). Let $d\in L(w)\setminus \{1\}$. Clearly, $C_1(w)=C_1(v)$. Since $C_d(w)$ and $C_1(w)$ touch at $w$, the cycle $C_d(w)$ contains the edges $uw$ and $ww_1$. Since $\phi(u)=2$, we have $d=2$ and the coloring is as shown in Figure \ref{fig:colors}(c$_2$).
\end{proof}

Let $Q_1,\dots,Q_{25}$ be the configurations shown in Figures \ref{fig:1}--\ref{fig:3.1}.
Our goal is to prove that each of these configurations is reducible.
We will use the notation about vertices of each of these configurations
as depicted in Figures \ref{fig:1}--\ref{fig:3.1} and in additional figures in this section.

\begin{lemma}
Configurations $Q_1$, $Q_2$, and $Q_3$ are reducible.
\end{lemma}

\begin{proof}
Let $\phi$ be an $L$-coloring of $D-v$. By Lemma
\ref{lem:triangle}(a), the cycle $C_1(v)$ in $Q_1$ uses the edges
$vv_2$ and $vv_3$, a contradiction to Lemma \ref{lem:directed
triangle}. Similarly, Lemma \ref{lem:triangle}(b) yields a
contradiction to Lemma \ref{lem:directed square} in the case of
$Q_2$. For $Q_3$, we apply Lemma \ref{lem:triangle}(b) at each of
the vertices $v,w,u$ and conclude that changing $\phi(w)$ to 2 and
$\phi(u)$ to 1 gives an $L$-coloring of $D-v$ that contradicts
Lemma \ref{lem:triangle}(a).
\end{proof}

\begin{lemma}
Configurations $Q_4$ and $Q_5$ are reducible.
\end{lemma}

\begin{proof}
We may assume that $vv_3 \in E(D)$. Let $\phi$ be an
$L$-coloring of $D-vv_3$. Clearly, we may assume that
$\phi(v)=\phi(v_3)=1$ for otherwise $\phi$ is an $L$-coloring of
$D$. The cycle $C=C_1(v)$ uses the arc $vv_3$. By Lemma
\ref{lem:directed triangle}, $C$ cannot use the arcs $v_2v$ or
$v_4v$. Therefore, we may assume that $C$ uses the arcs $v_1v$ and
$vv_3$. Since $C$ and the cycle $C'=C_2(v)$ touch at $v$, $C'$ uses the edges $v_4v$ and $vv_5$. In $Q_4$, this yields a contradiction by
Lemma \ref{lem:directed
triangle}. So, it remains to consider $Q_5$.

Let $u \neq v$ be the neighbor of $v_1$ on $C$ and $w \neq v$ be
the neighbor of $v_5$ on $C'$. Now, since $D$ is not
$L$-colorable and $\deg(v_1)=4$, if we were to recolor $v_1$ with color $c \in
L(v_1) \backslash \{1\}$ (and keep color 1 at $v$), we must have
$c=2$ and obtain a color-2 cycle $C''$ through $v_1$. Clearly,
$C''$ uses the edge $v_1v_5$. By Lemma \ref{lem:triangle}(a)
applied at the vertex $v_5$, $C''$ must use the edge $wv_5$.
Now, modifying $\phi$ by recoloring $v_5$ with the color $d\in L(v_5)\setminus \{2\}$ and $v$ with 2, we note that there cannot
be a color-2 cycle through $v$ nor a
color-$d$ cycle through $v_5$ since the only possibility for such a cycle is when $d=1$ and $C' \cup C''$ separates the
vertex $u$ from possible neighbor of $v_5$ of color 1 that is different from $v_1$.
Thus, we obtain an $L$-coloring of $D$, a contradiction.
\end{proof}

\begin{figure}[htb]
   \centering
   \includegraphics[width=8.6cm]{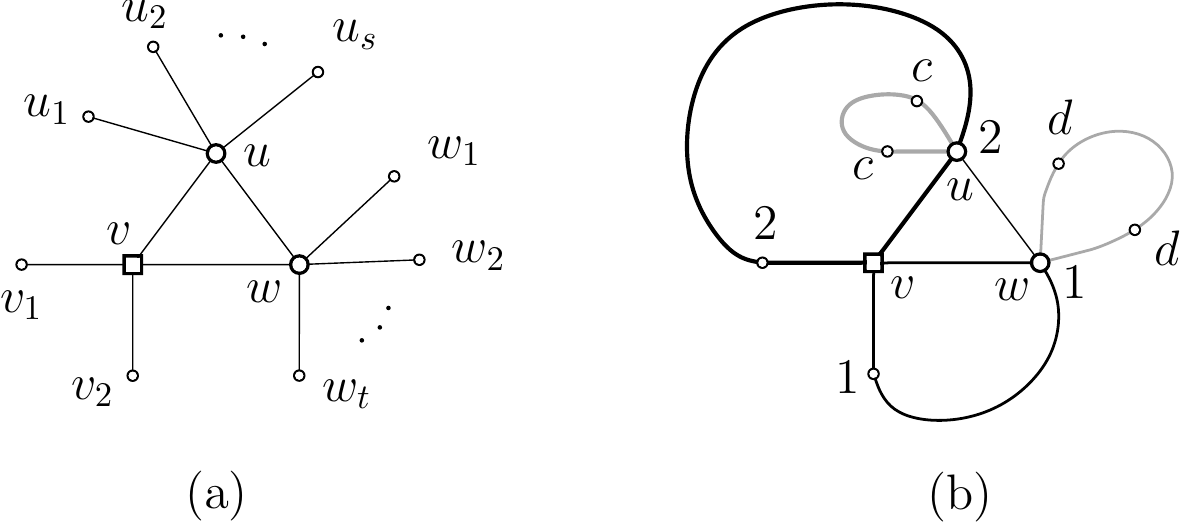}
   \caption{Triangle $T=vuw$ and its neighborhood}
   \label{fig:neighbors}
\end{figure}

In the proofs of all of the subsequent lemmas, showing reducibility of
particular configurations, we have a common scenario. Let us describe
the common notation and assumptions that we will inquire.

We will always have a triangle $T=vuw$, where $\deg(v)=4$. We shall
assume that $L(v)=\{1,2\}$ and will consider an $L$-coloring $\phi$ of $D-v$.
This coloring will also be denoted by $\phi_v$ if we would want to remind the
reader that the vertex $v$ is not colored.
The neighbors of the vertices of $T$ are denoted as in Figure
\ref{fig:neighbors}(a), $v_1,v_2$ being neighbors of $v$, $u_1,\dots,u_s$ neighbors of $u$ and $w_1,\dots,w_t$ neighbors of $w$, where $u_i$ and $w_j$ are enumerated in the clockwise order. It may be that $u_s=w_1$. By Lemma \ref{lem:triangle}(a), we may assume that
$\phi_v(v_2)=\phi_v(w)=1$ and $\phi_v(v_1)=\phi_v(u)=2$. We shall denote
the unused colors in $L(u)$ and $L(w)$ by $c$ and $d$, respectively, i.e.,
$c\in L(u)\setminus\{2\}$ and $d\in L(w)\setminus \{1\}$.
Sometimes we shall be able to conclude that $c=1$
or that $d=2$, but in general this needs not to be the case.

As discussed before, there are two cycles $C_1(v)$ and $C_2(v)$ passing
through $v$. Similar cycles can be defined for $u$ and $w$. First we
define an $L$-coloring $\phi_u$ of $D-u$ by modifying $\phi_v$ by
coloring $v$ with color 2 and uncolor $u$. This coloring defines cycles
$C_2(u)$ and $C_c(u)$ that touch at $u$.
Note that we may assume that $C_2(u)=C_2(v)$. Similarly, by coloring
$v$ with color 1 and uncolor $w$, we obtain a coloring $\phi_w$ of
$D-w$. The corresponding cycles $C_1(w)$ and $C_d(w)$ touch at $w$.
Note that $C_1(w)=C_1(v)$.
This situation is depicted in in Figure \ref{fig:neighbors}(b), where the touching of the cycles at $u$ and $w$ may be different than shown (e.g., the cycle $C_c(u)$ could be in the exterior of $C_2(u)$).
Note that if $c=1$ and $d=2$, it may happen that $C_c(u)$ and $C_d(w)$ share the edge $uw$ (but they would be disjoint elsewhere since $c\ne d$ in this case).

\begin{figure}[htb]
   \centering
   \includegraphics[width=14cm]{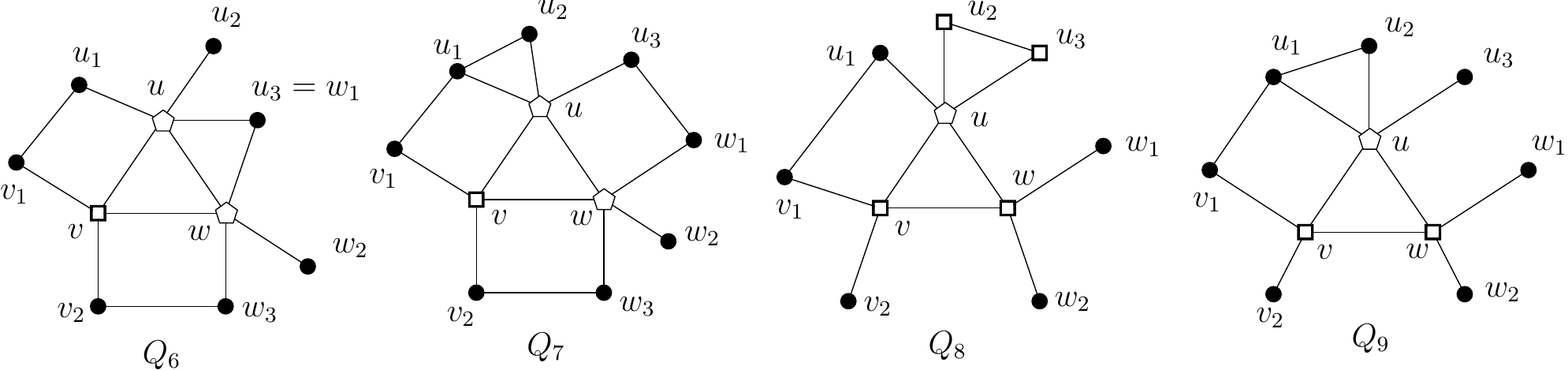}
   \caption{Configurations $Q_6$ to $Q_9$}
   \label{fig:Q6Q7Q8Q9}
\end{figure}

\begin{lemma}
Configurations $Q_6$, $Q_7$, $Q_8$, and $Q_9$ are reducible.
\end{lemma}

\begin{proof}
We shall use additional notation depicted in Figure \ref{fig:Q6Q7Q8Q9}.
We first assume that $C_1(v)$ uses the edge $ww_2$, and consider
the cycle $C_d(w)$. Since $C_1(v)$ and $C_d(w)$ touch at $w$,
$C_d(w)$ uses edges $wu$ and $ww_1$. In particular, this implies that $d=2$ since $\phi_w(u)=2$.
Note that this cannot happen in $Q_6$ by Lemma \ref{lem:directed triangle}, so we may assume to
have one of the configurations $Q_7$--$Q_9$.
Let us now consider the cycle $C_2(v)$. By Lemma \ref{lem:directed square}, $C_2(v)$ uses the
edge $u_2u$ or $u_3u$. If $C_2(v)$ uses the edge $u_2u$, the cycle $C_c(u)$ must use the edges $uw$ and $uu_3$. In particular, we have $c=1$. Now we recolor $u$ with color 1 and $w$ with color $d=2$. It is clear that there is no color-1 cycle through $u$ and there is no color-$d$ cycle through
$w$ (since it would need to touch $C_1(w)$). It follows that the modified coloring $\phi'$ is a proper $L$-coloring of
$D-v$ with the property that $\phi'(v_1)\ne\phi'(u)$, contradicting
Lemma \ref{lem:triangle}.

Thus, we may assume that $C_2(v)$ uses the edge $u_3u$, so $\phi_w(u_3)=\phi_v(u_3)=2$.
The cycle $C_d(w)$ uses the edges $uw$ and $ww_1$, and
hence we have $d=2$. Let us change
$\phi_v$ to a coloring $\phi'_v$ of $D-v$ by recoloring $u$ with color $c$ and $w$ with
color 2. A color-2 cycle through $w$ would touch $C_1(v)$ at the vertex $w$, and clearly,
there is no room for this, so there is no such cycle.
By Lemma \ref{lem:directed triangle}, there is no color-$c$ cycle through $u$ in $Q_7$
or $Q_9$. Therefore, the modified coloring $\phi'_v$ is a proper
$L$-coloring of $D-v$ for configurations $Q_7$ and $Q_9$. Now,
Lemma \ref{lem:triangle} yields a contradiction in these two
cases.

It remains to consider the configuration $Q_8$. If there were no
color-$c$ cycle through $u$, we would have a contradiction as
above. Therefore, $\phi(u_1) = \phi(u_2) = c$ and there is a
color-$c$ cycle $C'$ using the edges $u_1u$ and $uu_2$.
Now we modify the original coloring $\phi$ as follows. We uncolor $u_3$ and color $v$ with
color 2. Clearly, this gives an $L$-coloring $\phi'$ of $D-u_3$. 
By Lemma \ref{lem:triangle}(b), we conclude that $L(u_3)=L(u_2)=\{2,c\}$ and that $u_2$ has a neighbor of color 2 that is contained in the interior of the cycle $C'$. Now, we change $\phi'$ by coloring $u_3$ with color $c$ and $u_2$ with color 2.
By the above, it is easy to see that we obtain an $L$-coloring of $D$, a contradiction.  

Thus, it remains to consider the case that $C_1(v)$ does not use
the edge $w_2w$. By Lemma \ref{lem:triangle}(b), this case cannot
occur for the configurations $Q_8$ and $Q_9$. For the
configurations $Q_6$ and $Q_7$, $C_1(v)$ necessarily uses the edge
$w_1w$ by Lemma \ref{lem:triangle}(c). The cycle $C_d(w)$ must use
edges $ww_2$ and $ww_3$, thus
$\phi(w_2)=\phi(w_3)=d$. Now, consider $C_2(v)$. By Lemma
\ref{lem:directed square}, $C_2(v)$ cannot use the edge $u_1u$.
Assume that $C_2(v)$ uses the edge $u_2u$. Then $C_c(u)$ must use the
edges $u_3u$ and $uw$ in $Q_6$, and the edges $u_3u$, $uw$ and
$ww_1$ in $Q_7$. This contradicts Lemmas \ref{lem:directed
triangle} and \ref{lem:directed square}, respectively.
This settles the reducibility of $Q_6$.

For $Q_7$, assume finally that $C_2(v)$ uses the edge
$u_3u$. Now, $C_c(u)$ must use both edges $uu_1$ and $uu_2$,
which contradicts Lemma \ref{lem:triangle}. Hence, $Q_7$ is also reducible. This
completes the proof.
\end{proof}

\begin{figure}[htb]
   \centering
   \includegraphics[width=13cm]{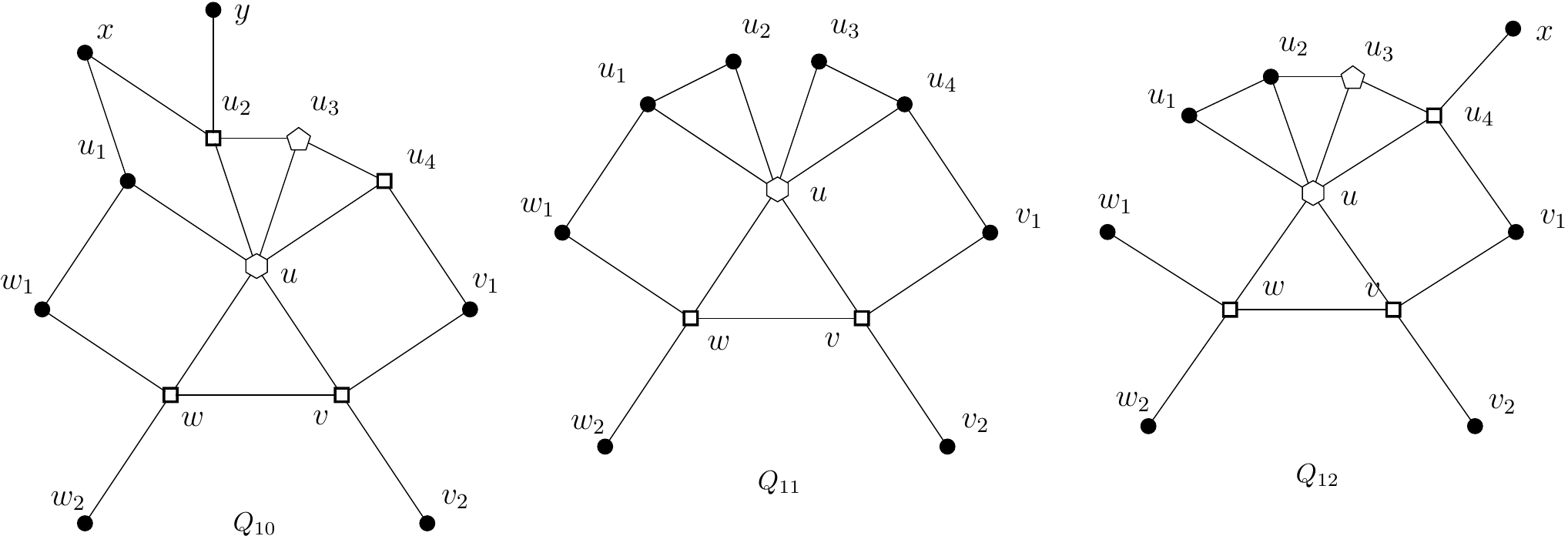}
   \caption{Configurations $Q_{10}$, $Q_{11}$, and $Q_{12}$}
   \label{fig:Q10Q11Q12}
\end{figure}

\begin{lemma}
Configurations $Q_{10}$, $Q_{11}$, and $Q_{12}$ are reducible.
\end{lemma}

\begin{proof}
We shall use additional notation depicted in Figure \ref{fig:Q10Q11Q12}.
By Lemma \ref{lem:triangle} we see that $L(w)=\{1,2\}$. Let $C=C_2(v)$. We may
assume that $C$ contains the directed arcs $uv$ and $vv_1$. The
cycle $C$ must use one of the arcs $u_1u$, $u_2u$, $u_3u$, or
$u_4u$. By Lemma \ref{lem:directed square}, $C$ cannot use the arc $u_4u$.
Next, suppose that $C$ uses the arc
$u_2u$. Now, we claim that modifying $\phi_v$ by recoloring $u$ with
color $c\in L(u)\setminus\{2\}$ and $w$ with color 2, gives an $L$-coloring of $D-v$.
Clearly, there is no color-2 cycle through $w$ since $w$ has only
one neighbor of color 2. Now, a color-$c$ cycle $C'$ through $u$ touches $C$, so it
uses the arcs $u_3u$ and $uu_4$, contradicting Lemma \ref{lem:directed triangle}.
Therefore, the modified coloring
is an $L$-coloring of $D-v$, and Lemma \ref{lem:triangle} now
implies that we can extend it to an $L$-coloring of $D$.

Now, suppose that $C$ uses the arc $u_1u$. Now, if we were to
modify $\phi_v$ by recoloring $u$ with color $c$ and $w$ with color 2,
by Lemma \ref{lem:triangle} this cannot be an $L$-coloring of $D-v$,
thus we must have a color-$c$ cycle $C'$ through $u$.
By Lemma \ref{lem:directed triangle}, $C'$ either uses the edges $u_2u$
and $u_4u$ or (in $Q_{11}$ only) uses the edges $u_2u$ and $uu_3$.
Suppose first that $C'$ uses the edges $u_2u$ and $uu_4$.
This cannot happen in $Q_{12}$ since it would contradict Lemma \ref{lem:triangle}(a)
at the vertex $u_4$ because in this case $C'$ would have to use the edge $u_4x$. 
Now, consider the cycle $C''=C_2(w)$ through $w$ in the coloring $\phi_w$. 
Then $C''$ must use the arcs $u_1u$, $uw$ and $ww_1$.
This contradicts Lemma \ref{lem:directed square} in cases $Q_{10}$ and $Q_{11}$.
The remaining case is that $C'$ uses
edges $u_2u$ and $uu_3$, which can happen only in $Q_{11}$
(by Lemma \ref{lem:directed square}). The cycle $C_2(w)$, which uses edges
$w_1w$ and $wu$, cannot use $u_1u$ by Lemma \ref{lem:directed square}, so it must
use the edge $uu_4$.

Next, we distinguish two cases. First, assume that $u_4u, uw, ww_1 \in
E(D)$. Then $u_4v_1 \in E(D)$. Clearly, $C_2(w)$ uses a vertex $x\ne u$ on
$C$ (since $C$ and $C_2(w)$ cross at $u$). But now, the directed path from $v_1$ to $x$ on $C$, together with the
arc $u_4v_1$ and the path from $x$ to $u_4$ on $C_2(w)$ create a
directed color-2 closed walk in the original coloring $\phi$, a
contradiction. Secondly, assume that $uu_4, wu, w_1w \in E(D)$.
Clearly, $C_2(w)$ uses a vertex $y\ne u$ on $C$. But now, the directed
path from $y$ to $u_1$ on $C$, with the arcs $u_1u, uu_4$ and the
directed path from $u_4$ to $y$ on $C_2(w)$ creates a color-2
directed closed walk in the original coloring $\phi$, a contradiction.

It remains to consider the case when $C$ uses the arc $u_3u$. Now, if we were to
modify $\phi$ by recoloring $u$ with color $c$ and $w$ with color
2, then Lemma \ref{lem:triangle} would imply that this is not an
$L$-coloring of $D-v$. Since there cannot be a color-2 cycle through
$w$, this implies that there is a color-$c$ cycle $C'$ through $u$.
Since $C$ and $C'$ touch, $C'$ must use the
edges $u_2u$ and $uu_1$, and hence $\phi(u_1)=\phi(u_2)=c$.
By Lemma \ref{lem:directed triangle}, this is not possible in $Q_{11}$ and
$Q_{12}$, so we are in $Q_{10}$. Since
$C'$ is a directed cycle, assume that we have the arcs $u_2u$ and
$uu_1$ (similar argument works for the other possibility). Now,
$C'$ cannot use the arc $xu_2$ by Lemma \ref{lem:directed square}.
Therefore, $C'$ uses the arc $yu_2$. Now, modify $\phi$ by
recoloring $u$ with color $c$, $w$ with color 2, color $v$ with color 1, and uncolor $u_2$.
The resulting coloring of $D-u_2$ contradicts Lemma \ref{lem:triangle}(a).
This completes the proof.
\end{proof}

\begin{figure}[htb]
   \centering
   \includegraphics[width=13.3cm]{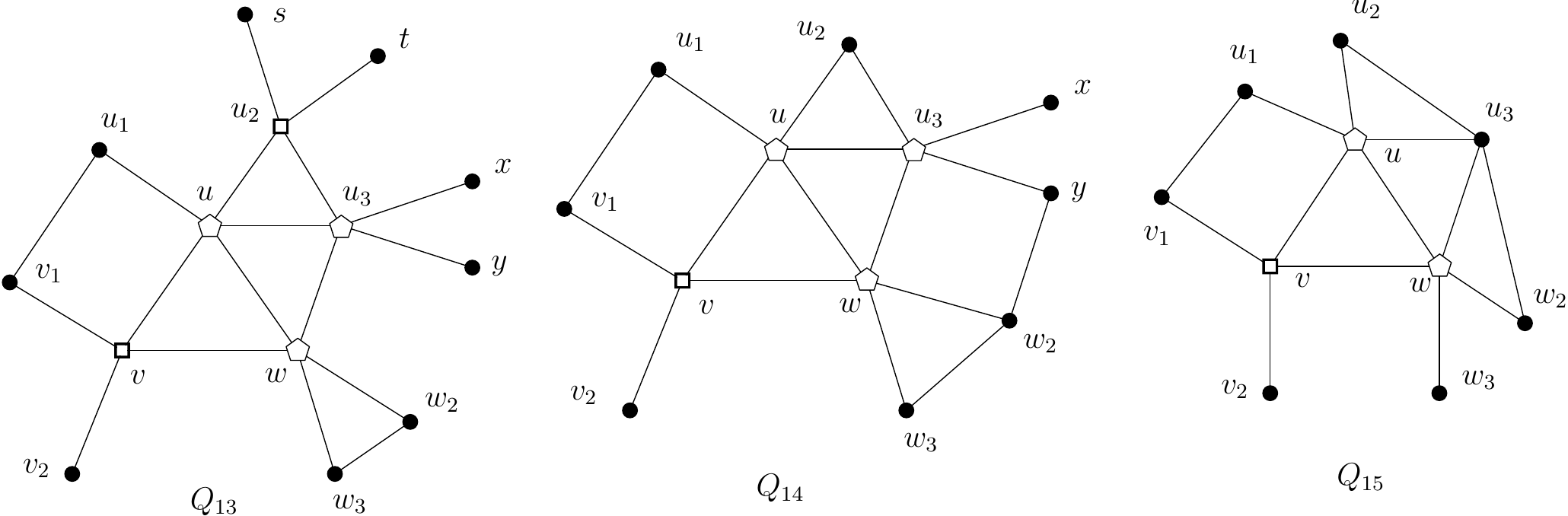}
   \caption{Configurations $Q_{13}$, $Q_{14}$, and $Q_{15}$}
   \label{fig:Q13Q14Q15}
\end{figure}

\begin{lemma}
Configurations $Q_{13}$, $Q_{14}$ and $Q_{15}$ are reducible.
\end{lemma}

\begin{proof}
The cycle $C_2(v)$ uses the edges $uv$ and
$vv_1$ and we may assume $uv, vv_1 \in E(D)$. By Lemma
\ref{lem:directed square}, $C_2(v)$ cannot use the edge $u_1u$.
If $C_2(v)$ uses the edge $u_2u$, then $C_c(u)$ uses the edges $u_3u$ and
$uw$, contradicting Lemma \ref{lem:directed triangle}.
Therefore, we may assume henceforth that $C_2(v)$ uses the edge $u_3u$
and that, in particular, $\phi(u_3)=2$. Therefore, the cycle $C_c(u)$ uses the edges
$u_1u$ and $uu_2$, and we have that $\phi_v(u_1) = \phi_v(u_2)=c$.

The cycle $C''=C_d(w)$ uses two of the incident arcs to $w$.
Since two neighbors of $w$ are on $C_1(v)$ and $u,u_3$
are on $C_2(v)$, we conclude that $d=2$.
Clearly, by Lemma \ref{lem:directed triangle}, $C''$ cannot use
both of the edges $wu_3$ and $uw$. Since $C_1(v)$ and $C''$ touch at $w$,
$C_1(v)$ contains $w_3$ and $C''$ contains the edge $ww_2$.

Note that $u_3\in V(C'')$, since $u_3$ is the
only neighbor of $u$ of color 2 in the coloring $\phi_w$.
Let us first suppose that $C''$
contains the arc $uw$. Since $D$ has no directed triangles, we conclude that
$u_3w\in E(D)$, so we may shorten $C''$ by eliminating vertex $u$, and thus
we may henceforth assume that $C''$ contains the edge $u_3w$.
Now, Lemma \ref{lem:directed triangle} yields a contradiction in the case of the
configuration $Q_{15}$. So, we are left to consider $Q_{13}$ and $Q_{14}$.

Suppose that $\phi_v(x)=2$, where $x$ is the neighbor of $u_3$ as shown in the figure.
Then we modify the original coloring $\phi$ as follows: we recolor
$u_3$ with the color $c'\in L(u_3)\setminus \{2\}$, recolor $w$ with color 2, and color $v$ with
color 1. We claim that this is an $L$-coloring of $D$. Clearly, there
is no color-1 cycle through $v$ since $v_2$ is the only neighbor
of $v$ with color 1. Similarly, there is no color-2 cycle through
$w$ since $u$ has no neighbor of color 2.
Lastly, there is no color-$c'$ cycle through $u_3$, since such a cycle would
need to use the edges $u_3u_2$ and $u_3y$, thus it would not touch $C''$.
This contradiction shows that $C_2(v)$ and $C''$ use the edge $yu_3$, and consequently,
$\phi(y)=2$. In particular, the cycle $C''=C_2(w)$ uses the edges $yu_3,u_3w$, and
$ww_2$. Lemma \ref{lem:directed square} yields contradiction in the case of the
configuration $Q_{14}$.

It remains to consider $Q_{13}$.
First, we observe that the color-$c$ cycle $C'=C_c(u)$ mentioned above uses the edge $u_2s$.
(This follows by Lemma \ref{lem:triangle}(a) applied to the coloring of $D-u_2$
obtained in this case.) Therefore, $\phi(s)=c$ and $\phi(t)=\phi(u_3)=2$.

Now, if we were to modify the original coloring $\phi$ by
recoloring $u_3$ with the color $c'$, we would obtain a
color-$c'$ cycle $Q$ through $u_3$ for otherwise coloring
$v$ with color 2 would give a proper $L$-coloring of $D$.
Clearly, $Q$ and $C_2(v)$ touch at $u_3$, so $Q$ uses edges $u_2u_3$ and $u_3x$.
In particular, we have $c=c'$. But now we can recolor $u_2$ with the
color $c''\in L(u_2)\setminus \{c\}$, $v$ with color 2 and obtain an $L$-coloring of $D$.
This contradiction completes the proof.
\end{proof}

\begin{figure}[htb]
   \centering
   \includegraphics[width=13.5cm]{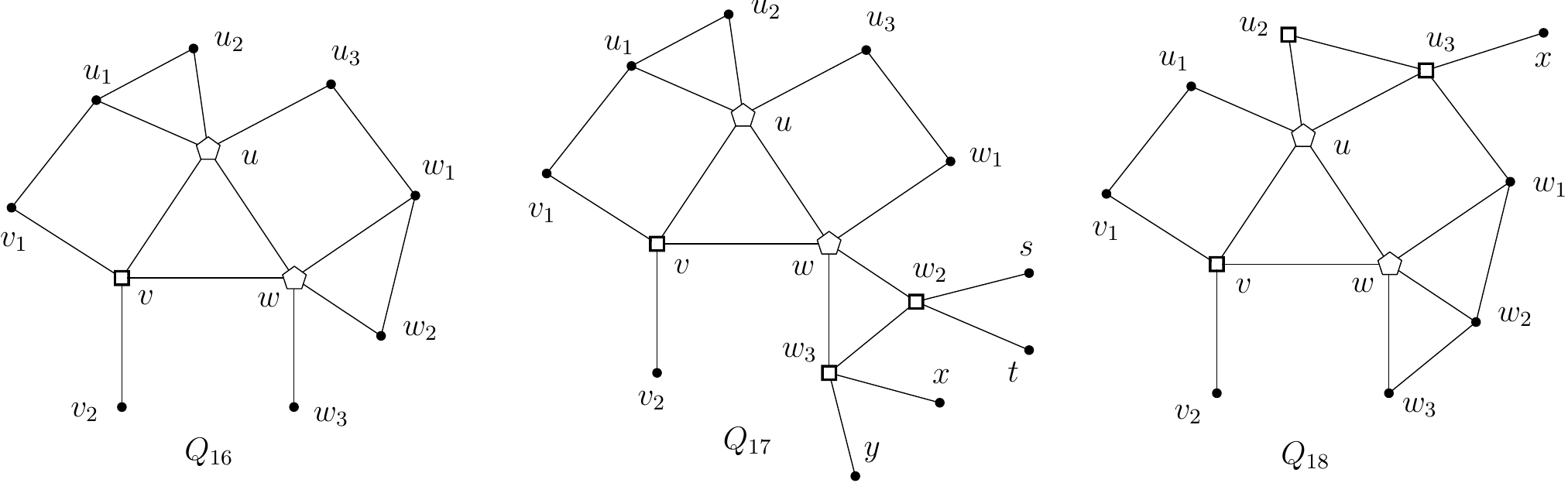}
   \caption{Configurations $Q_{16}$, $Q_{17}$, and $Q_{18}$}
   \label{fig:Q16Q17Q18}
\end{figure}

\begin{lemma}
The configurations $Q_{16}, Q_{17}, Q_{18}$ are reducible.
\end{lemma}

\begin{proof}
We may assume $uv, vv_1 \in E(D)$.
By Lemma \ref{lem:directed square}, $C_2(v)$ cannot use the edge $u_1u$.
Therefore, $C_2(v)$ uses one of the edges $u_2u$ and $u_3u$.

Let us first assume that $C_2(v)$ uses the edge $u_3u$.
Then the cycle $C_c(u)$ uses edges $u_1u$ and $uu_2$.
Lemma \ref{lem:directed triangle} gives a contradiction in the case of
configurations $Q_{16}$ and $Q_{17}$, so it remains to consider $Q_{18}$.
Let us now modify $\phi_v$ as follows: color $v$ with
color 2 and uncolor $u_3$. Clearly, this is a proper $L$-coloring of
$D-u_3$. Thus, Lemma \ref{lem:triangle}(a) implies that
$\phi_v(w_1)=\phi_v(u)=2$, and that $\phi_v(x)=\phi_v(u_2)=c$.

Now, we consider the cycle $C''=C_2(w)$.
Since $C''$ and $C_1(w)$ touch at $w$,
$C''$ cannot use the edge $ww_3$. By Lemma \ref{lem:directed triangle},
$C''$ cannot use both edges $w_1w$ and $ww_2$.
Similarly, $C''$ cannot use both, $uw$ and $ww_1$, since in that case $C''$
would also use $u_3u$, contradicting Lemma \ref{lem:directed square}. 
The only possibility left is that $C''$ uses the edges $uw$ and $ww_2$. 
This implies that
$C''$ uses the arcs $w_1u_3$, $u_3u$, $uw$, and $ww_2$. Consequently,
$w_1w \in E(D)$. Since $D$ has no
directed triangles, it follows that $w_1w_2 \in E(D)$.
Let $P$ be the directed path on $C''$ from $w_2$ to $w_1$.
But now the directed closed walk $w_1w_2P$ contains a color-2
cycle in the original coloring $\phi_v$, a contradiction. This completes the proof
when $u_3u$ belongs to $C_2(v)$.

Suppose now that $C_2(v)$ uses the edge $u_2u$.
Then the cycle $C'=C_c(u)$ uses
the edges $u_3u$ and $uw$. Note that this implies in particular that
$c=1$. By Lemma \ref{lem:directed square},
$C'$ cannot use the edge $ww_1$. Thus, we have two possibilities:
$C'$ uses the edge $ww_2$ or $ww_3$.

In either case, we modify the original coloring $\phi_v$ by recoloring
$u$ with color $c=1$, coloring $v$ with color 2, and uncoloring $w$.
We now consider the cycles $C_1(w)$ and $C_d(w)$ with respect to this coloring of $D-w$.
Clearly, $C'=C_1(w)$. The cycle $C''=C_d(w)$ touches $C'$ at $w$.
If $C'$ used the edge $ww_2$, then $C_d(w)$ would have to use edges
$ww_3$, $wv$ and $vv_1$, which is not possible since it would need to cross
the color-1 cycle $C_1(v)$. Thus, $C'$ uses the edge $ww_3$, and $C''$
uses the edges $w_1w$ and $ww_2$. In cases $Q_{16}$ and $Q_{18}$
we have a contradiction to Lemma \ref{lem:directed triangle}.
This completes the proof for $Q_{16}$ and $Q_{18}$.

It remains to consider $Q_{17}$.
Recall that $C'$ uses the edge $ww_3$. Note
that $C'$ cannot use the edge $w_3x$ since then we could modify
the original coloring $\phi$ by recoloring $u$ with color 1, $w_3$
with color $d'\in L(w_3)\setminus \{1\}$, and coloring $v$ with color 2, obtaining a proper
$L$-coloring of $D$. Therefore, $C'$ uses the edge $w_3y$ and,
consequently, $\phi(w_3)=\phi(y)=1$.
Since $w_1,w_2$ belong to $C''=C_d(w)$, we have
$\phi(w_1) = \phi(w_2)=d$. Now, $C''$ uses the edge
$w_2s$ (by Lemma \ref{lem:triangle}(b)). In particular, we have
$\phi(s)=d$, and consequently, $\phi(t)=1$.

Now, we modify the original coloring $\phi$ by recoloring $u$ with
color 1 and $w_3$ with color $d'$. Since this is not a proper
$L$-coloring of $D-v$, it follows that there is a color-$d'$ cycle
$C'''$ through $w_3$. Clearly, $C'''$ must use the edges
$sw_2,w_2w_3$ and $w_3x$, separating $t$ from $w$. Now, recoloring
$w_2$ with color 1 and coloring $v$ with color 2, we obtain a
proper $L$-coloring of $D$. This final contradiction shows that
$Q_{17}$ is reducible.
\end{proof}

\begin{figure}[htb]
   \centering
   \includegraphics[width=10cm]{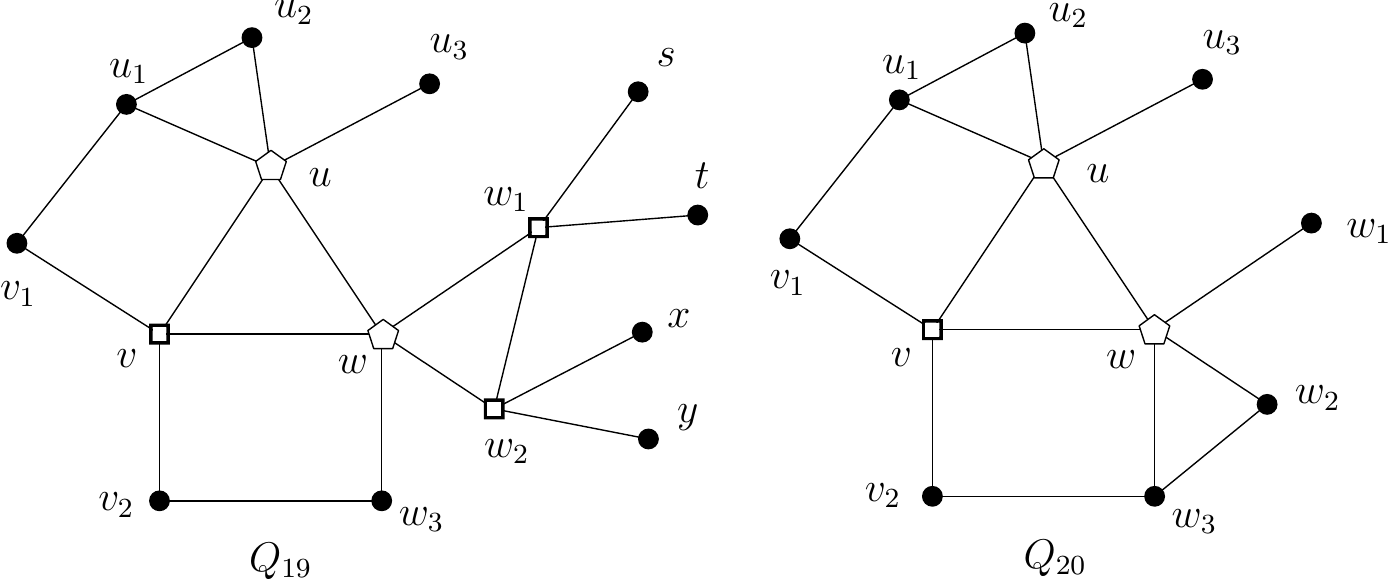}
   \caption{Configurations $Q_{19}$ and $Q_{20}$}
   \label{fig:Q19Q20}
\end{figure}

\begin{lemma}
The configurations $Q_{19}$ and $Q_{20}$ are reducible.
\end{lemma}

\begin{proof}
We may assume that $uv, vv_1 \in E(D)$.
The cycle $C_2(v)$ cannot use the edge $u_1u$ by Lemma
\ref{lem:directed square}. Therefore, $C_2(v)$ must use one of the edges
$u_2u$ or $u_3u$. If $C_2(v)$ uses the edge $u_3u$, then $C_c(u)$ would use
edges $u_1u$ and $uu_2$, a contradiction to Lemma \ref{lem:directed triangle}.
Thus, $C_2(v)$ uses the edge $u_2u$. The cycles $C_c(u)$ and $C_2(v)$ touch at $u$; thus $C_c(u)$
must use the edges $uw$ and $uu_3$. This implies that $c=1$.
Let us now consider the cycle $C_1(v)$. Clearly, it contains the edges
$wv$ and $vv_2$, and does not contain the
edge $w_3w$ by Lemma \ref{lem:directed square}. Therefore, $C_1(v)$
uses one of the edges, $w_1w$ or $w_2w$.

If $C_1(v)$ uses the edge $w_2w$, then $C_d(w)$ must use the edges $wu$ and $ww_1$.
Since $\phi_w(u)=2$, we have $d=2$ and $\phi(w_1)=2$. Observe that the cycles
$C_1(u)$ and $C_2(w)$ share the edge $uw$, but are otherwise disjoint. However, this 
is not possible, since they ``cross each other'' when viewed how they leave the edge
$uw$ at one and the other end. This contradiction shows that $C_1(v)$ uses the edge $w_1w$.
Consequently, we have $\phi(w_1)=1$. Then $C_d(w)$ contains the edges
$ww_3$ and $ww_2$. This contradicts Lemma \ref{lem:directed triangle}
for configuration $Q_{20}$.

It remains to consider $Q_{19}$. As mentioned above,
$C_1(v)$ uses the edge $w_1w$. If $C_1(v)$ uses the edge $tw_1$,
we modify the original coloring $\phi$ by coloring $v$ with color
1 and uncoloring $w_1$. Clearly, this modified coloring is a
proper $L$-coloring of $D-w_1$. But now, Lemma \ref{lem:triangle}(a)
yields a contradiction.
Hence, we may assume that $C_1(v)$ uses the edge $sw_1$ and consequently $\phi(s)=1$. Now, if were
to modify the original coloring $\phi$ by coloring $v$ with color
1, recoloring $w$ with color $d$, and uncoloring $w_2$, we would
obtain a proper $L$-coloring of $D-w_2$. By Lemma \ref{lem:triangle}(a), it follows that
$L(w_2)=\{1,d\}$ and that $\phi_v(y)=d$ and $\phi_v(x)=1$.

Now, we modify the original coloring $\phi_v$ as follows: we color
$v$ with color 1, recolor $w$ with color $d$, recolor $w_2$ with
color 1 and uncolor $w_1$. Since $C_1(v)$ separates $u$ from $w_3$
there is no color-$d$ cycle through $w$. It is easy to
see that there are no monochromatic cycles through $v$ or $w_2$.
Therefore, the modified coloring is a proper $L$-coloring of
$D-w_1$. But now, Lemma \ref{lem:triangle} implies that
$\phi(s)=d$, a contradiction. Therefore, $Q_{19}$ is also reducible.
\end{proof}

\begin{figure}[htb]
   \centering
   \includegraphics[width=13.3cm]{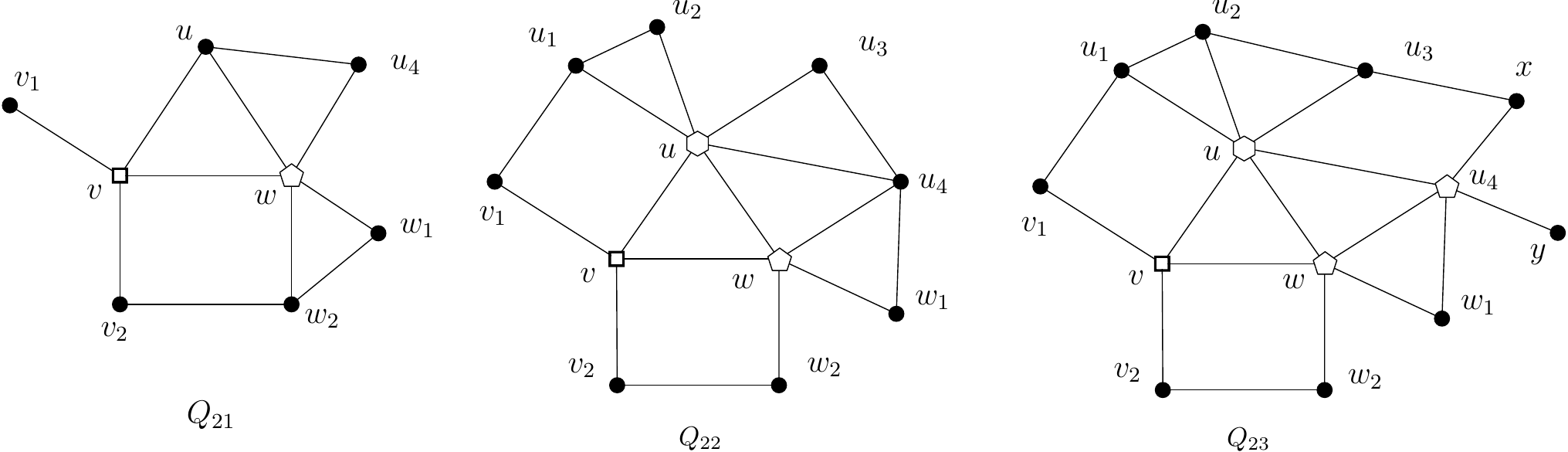}
   \caption{Configurations $Q_{21}$, $Q_{22}$, and $Q_{23}$}
   \label{fig:Q21Q22Q23}
\end{figure}

\begin{lemma}
The configurations $Q_{21}$, $Q_{22}$ and $Q_{23}$ are reducible.
\end{lemma}

\begin{proof}
Consider the cycle $C_1(v)$. It uses the edges $wv$ and $vv_2$.
By Lemma \ref{lem:directed square}, $C_1(v)$ cannot use the edge $w_2w$.
If $C_1(v)$ uses the arc $w_1w$, then $C_d(w)$ would use edges
$uw$ and $wu_4$, and Lemma \ref{lem:directed triangle} would yield
a contradiction. This proves that $C_1(v)$ uses the arc $u_4w$.
Now we see that the cycle $C_d(w)$ uses the edges
$w_1w$ and $ww_2$, and thus $\phi_v(w_1) = \phi_v(w_2) = d$.
Lemma \ref{lem:directed triangle} gives a contradiction in the case of $Q_{21}$,
so it remains to consider $Q_{22}$ and $Q_{23}$.

The cycle $C_2(v)$ uses the edges $uv$ and
$vv_1$. By Lemma \ref{lem:directed square}, it cannot use the
edge $u_1u$. Since $\phi_v(u_4)=1$, $C_2(v)$ uses one of the edges
$u_2u$ or $u_3u$. If $C_2(v)$ uses the edge $u_3u$, consider the cycle $C_c(u)$.
It cannot use the edges $u_1u$ and $uu_2$, since this would
contradict Lemma \ref{lem:directed triangle}. Similarly, it cannot use the
edges $uw$ and $uu_4$. This shows that $C_2(v)$ cannot use the edge $u_3u$.

We conclude from the above that $C_2(v)$ uses the edge $u_2u$.
The cycles $C_2(u)=C_2(v)$ and $C'=C_c(u)$ touch at $u$. Therefore,
$C'$ uses two of the edges $uu_3,uu_4,uw$, and in particular we have $c=1$
since $\phi_u(w)=\phi_u(u_4)=1$.

In $Q_{22}$, $C'$ cannot use both of the edges $u_3u$ and $uu_4$ by
Lemma \ref{lem:directed triangle}. Similarly, $C'$ cannot use
both of the edges $u_4u$ and $uw$. Therefore, we may assume that
$C'$ uses the edges $u_3u$ and $uw$. But then, $C'$ necessarily uses
the edges $wu_4$ since $\phi_u(w_1)=\phi_u(w_2)=d$ and $\phi_u(v)=2$.
Since we may assume that $u_3uwu_4$ is a directed path, we get a
contradiction with Lemma \ref{lem:directed square}.

It remains to consider $Q_{23}$.
Recall that $C_2(v)$ uses the edge $u_2u$, and also
note that $C_1(v)$ uses one of the
arcs $yu_4$ or $xu_4$. First, suppose that $C_1(v)$ uses the arc $yu_4$ so that
$\phi_v(y)=1$. Now, we modify the coloring $\phi_u$ by
recoloring $u_4$ with color $d'\in L(u_4)\setminus \{1\}$
and coloring $u$ with color 1. Now, since $C_1(v)$
separates $x$ from $w_1$, there is no color-$d'$ cycle through $u_4$
in the modified coloring. Similarly, since $C_2(v)$ separates $u_1$
from $u_3$ and since $\phi_u(w_1)=\phi_u(w_2)=d$, there is no color-1 cycle
through $u$ in the modified coloring. This shows that we have a proper
$L$-coloring of $D$, a contradiction.

Thus, we may assume that $C_1(v)$ uses the arc $xu_4$, and
consequently, $\phi_v(x)=1$. Now, we claim that $\phi_v(y)=d$. If not,
then by the argument above we could modify $\phi_u$ by recoloring
$u_4$ with color $d'$ and coloring $u$ with color 1, obtaining a proper
$L$-coloring of $D$ unless $d'=d$ and there is a color-$d$ cycle
using edges $yu_4$ and $u_4w_1$. This proves our claim that $\phi_v(y)=d$.

Now, we again consider the cycle $C' = C_c(u)$ (recalling that $c=1$ and that
$C'$ uses two of the three edges $u_3u$, $u_4u$ and $wu$ at
the vertex $u$. If $C'$ uses the edges $u_3u$ and $uu_4$, then since
$\phi_v(w_1)=\phi_v(w_2)= \phi_v(y)=d$, $C'$ must also
use the edge $u_4x$, contradicting Lemma \ref{lem:directed
square}. Clearly, by Lemma \ref{lem:directed triangle}, $C'$
cannot use both of the edges $u_4u$ and $uw$. Therefore, it
remains to consider the case that $C'$ uses the edges $u_3u$ and
$uw$. But then, $C'$ necessarily uses the edges $wu_4$ and
$u_4x$. Since $D$ has no directed triangles, it follows that
$u_3uu_4x$ is a directed path of color 1, implying that there was a
color-1 cycle through the edge $u_3x$ in the original coloring
$\phi$, a contradiction. This completes the proof of reducibility of $Q_{23}$.
\end{proof}

\begin{figure}[htb]
   \centering
   \includegraphics[width=11.2cm]{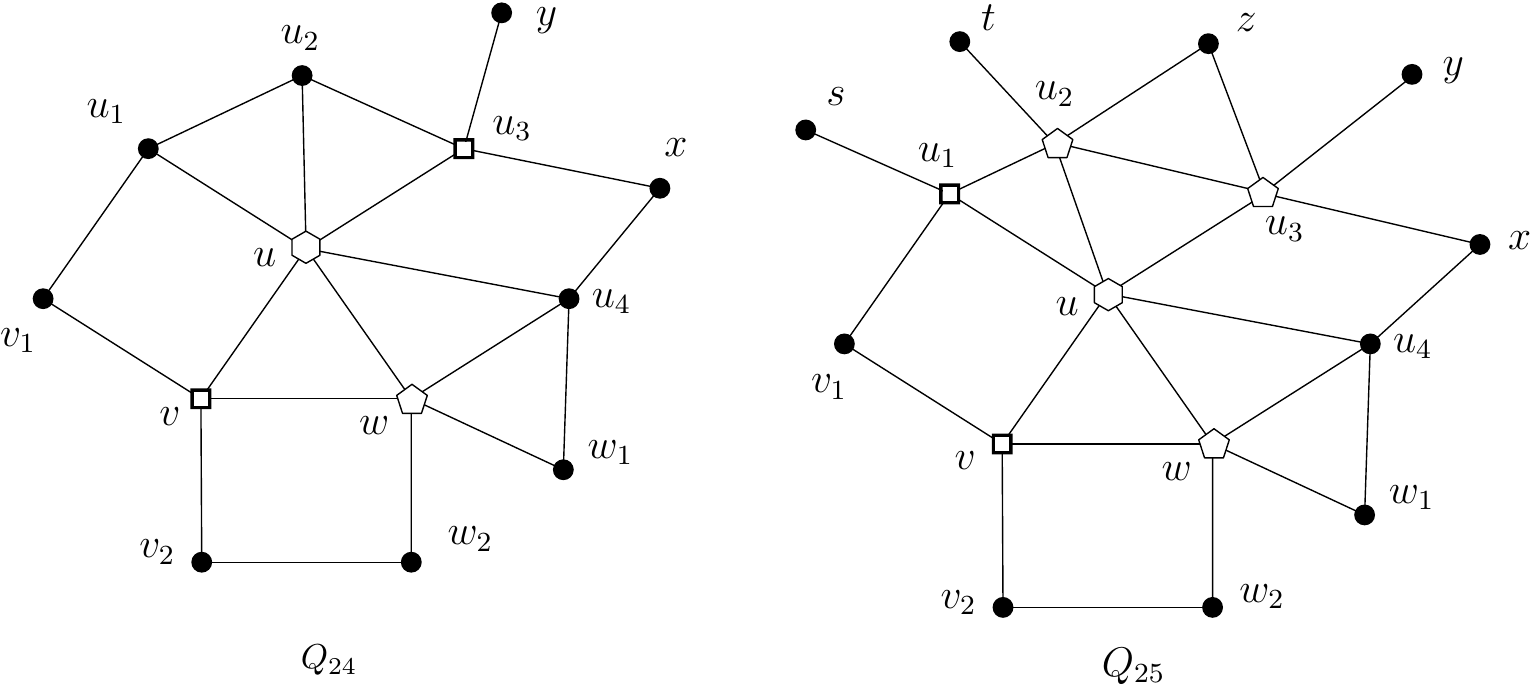}
   \caption{Configurations $Q_{24}$ and $Q_{25}$}
   \label{fig:Q24Q25}
\end{figure}

\begin{lemma}
The configurations $Q_{24}$ and $Q_{25}$ are reducible.
\end{lemma}

\begin{proof}
The cycle $C=C_1(v)$ uses the edges $wv$ and $vv_2$.
By Lemma \ref{lem:directed square}, $C$ cannot use the edge $w_2w$.
If $C$ uses the edge $ww_1$, then the cycle $C_d(w)$, which touches $C$
at $w$, has to use the edges $uw$ and $wu_4$. However, this contradicts
Lemma \ref{lem:directed triangle}.
Therefore, $C$ uses the arc $u_4w$ and $\phi(u_4)=1$. Now, consider $\phi_w$.
Clearly, $C_d(w)$ uses the edges $ww_1$ and $ww_2$.
It follows that $\phi_v(w_1)=\phi_v(w_2)=d$. Hence, since $d \neq 1$, $C$ cannot use the edge $u_4w_1$.

Now, consider the cycle $C'=C_2(v)$.
Since $C'$ uses the edges $uv$ and $vv_1$, it cannot use the edge $u_1u$
by Lemma \ref{lem:directed square}. Since $\phi(u_4)=1$, $C'$ must use one of the edges
$u_2u$ or $u_3u$.
If $C'$ uses the edge $u_3u$, consider the coloring $\phi_u$ and the cycle $C_c(u)$.
This cycle touches $C_2(u)$ at $u$. By Lemma \ref{lem:directed triangle},
it can neither use the edges $uu_1$ and $uu_2$ nor the edges $u_4u$ and $uw$.
This contradiction proves that $C'$ uses the edge $u_2u$.
The cycle $C_c(u)$, which touches $C'$ at $u$, must use the edge $u_3u$.
Since $\phi_u(w)=\phi_u(u_4)=1$, we have $c=1$.
Also note that $\phi_u(v)=2$; thus, if $C_c(u)$ uses the edge $uw$, then it also uses the edge $wu_4$,
and we can replace the cycle by a shorter color-1 cycle by using the edge $uu_4$ instead of $uw$ and $wu_4$.
(The orientation of the path $uwu_4$ is the same as $uu_4$ since $D$ has no directed 3-cycles.)
Thus, we may assume that $C_c(u)$ contains the edge $uu_4$. By Lemma \ref{lem:directed square},
$C_c(u)$ cannot use the edge $u_3x$, so it must use the edge $u_3y$ (or $u_3z$).

Let us now consider $Q_{24}$. Changing $\phi_u$ by uncoloring $u_3$ and coloring $u$ with color 1,
we obtain an $L$-coloring of $D-u_3$. This coloring is in contradiction with
Lemma \ref{lem:triangle}(a) at the vertex $u_3$, so this concludes the proof for $Q_{24}$.

It remains to show reducibility of $Q_{25}$. Let us recall that $C'$ uses the edge $u_2u$.
If $C'$ used the edge $u_2u_1$, then we could replace its edges $uu_2$ and $u_2u_1$ by the
edge $uu_1$, and this would lead to a contradiction to Lemma \ref{lem:directed square}.
Clearly, $C'$ cannot use the edge $u_2u_3$ (since $\phi_u(u_3)=1$), so it
either uses the edge $u_2t$ or $u_2z$.

First, assume that $C'$ uses the edge $u_2t$.
Now, we modify the coloring $\phi_v$ by recoloring $u_2$
with color $d' \in L(u_2) \backslash \{2\}$, and coloring $v$ with color 2. Since we previously showed
that any color-2 cycle through $v$ must use the edge $u_2u$, it
follows that there is no color-2 cycle through $v$. Since $D$ is
not $L$-colorable, it follows that there is a color-$d'$ cycle $C''=C_{d'}(u_2)$ through
$u_2$. Since $C'$ and $C''$ touch at $u_2$, it follows
that $C''$ contains the edges $zu_2$ and $u_2u_3$,
contradicting Lemma \ref{lem:directed triangle}.

Finally, assume that $C'$ uses the edge $zu_2$.
Now, extend the original coloring $\phi$ by coloring $v$ with
color 2 and recoloring $u_2$ by the color $d' \in L(u_2) \backslash \{2\}$.
This yields a color-$d'$ cycle $C_{d'}(u_2)$ which touches $C'$ at $u_2$ and thus uses the edges
$tu_2,u_2u_1$, and $u_1s$. It follows
that $\phi(u_1)= \phi(t) = \phi(s)=d'$. Now, recoloring $u_1$ with
color $c' \in L(u_1) \backslash \{d'\}$, we get a proper $L$-coloring of $D$
(by using Lemma \ref{lem:directed square}).
This final contradiction shows that $Q_{25}$ is reducible.
\end{proof}

We are ready to complete the proof of the main result.

\begin{proof}[Proof of Theorem \ref{thm:main choosable}]
By Theorem \ref{thm:unavoidable}, every planar graph of
minimum degree at least four contains one of the
configurations $Q_1,\dots,Q_{25}$. Suppose that $D$ is a minimum
counterexample to Theorem \ref{thm:main choosable}. Then $D$ has
digirth at least five and minimum degree at least four, but cannot
contain any of the configurations $Q_1,\dots,Q_{25}$ by Lemmas
3.5-3.13. This proves that a counterexample
does not exist, and the proof is complete.
\end{proof}

\section{Concluding remarks}

We raise the following questions. It would be interesting to see if
the result can be pushed to digirth 4. Also, the following relaxation of Conjecture \ref{conj:planar}
should be of interest.

\begin{conjecture} There exists $k$ such that every simple planar digraph
without cycles of length $4,\dots,k$ is 2-colorable.
\end{conjecture}

The original conjecture still seems out of reach. In fact, we do not know of a simple proof of the fact that
planar digraphs of large digirth are 2-colorable. In support of the conjecture, it would be nice to see whether one can find large acyclic set in a planar digraph, say of size $n/2$.
In fact, the following was conjectured in \cite{H2011}.

\begin{conjecture}
\label{conj:acyclic set}
Every simple $n$-vertex planar digraph has an acyclic set of size at least
$\tfrac{3n}{5}$.
\end{conjecture}

It is known that the bound in Conjecture \ref{conj:acyclic set} cannot be replaced by any larger value.

\end{document}